\documentclass[12pt,a4paper,english]{article}
\usepackage[T1]{fontenc}
\usepackage[latin1]{inputenc}
\usepackage{graphicx}

\makeatletter

\usepackage{amsthm,amsfonts,amsmath,amssymb}

\theoremstyle{definition}
\newtheorem{defn}{Definition}
\theoremstyle{plain}
\newtheorem{thm}[defn]{Theorem}
\theoremstyle{plain}
\newtheorem{lem}[defn]{Lemma}
\theoremstyle{plain}
\newtheorem{cor}[defn]{Corollary}
\theoremstyle{definition}
\newtheorem*{rem*}{Remark}
\theoremstyle{plain}
\newtheorem{observation}[defn]{Observation}

\usepackage{fullpage,enumerate}

\newcommand{\dtv}{\mathrm{d}_\mathrm{TV}}

\def\epsilon{\varepsilon}

\def\Ham{\textnormal{Ham}}

\newcommand{\E}[1]{\textnormal{\textbf{E}} \left[ #1 \right]}

\def\Prob{\textnormal{Pr}}

\newcommand{\Z}{\mathbb{Z}}
\newcommand{\1}{\mathbf{1}}

\def\Pk{P^{[k]}}

\def\P1{P^{[1]}}
\def\Pm{P^{[m]}}
\def\Ham{\textnormal{Ham}}

\def\M{\mathcal{M}}
\def\scan{\mathcal{M}_\rightarrow}

\def\HscanDob{\mathcal{M}_\textnormal{AnyOrder}}
\def\HscanPC{\mathcal{M}_\textnormal{FixedOrder}}
\def\scan{\mathcal{M}_\rightarrow}
\def\Hrandom{\mathcal{M}_\textnormal{RND}}

\def\M{\mathcal{M}}
\def\Mix{\textnormal{Mix}}

\def\HscanPCP{P_\textnormal{FixedOrder}}

\def\state{\Omega_\sim}

\usepackage{babel}
\makeatother
\begin{document}

\title{On systematic scan for sampling $H$-colourings of the path%
\thanks{This work was partly funded by EPSRC project GR/S76168/01. %
}}

\author{Kasper Pedersen\\
Department of Computer Science\\
University of Liverpool\\
Liverpool L69 3BX, UK\\
k.pedersen@csc.liv.ac.uk}
\maketitle
\begin{abstract}
This paper is concerned with sampling from the uniform distribution
on $H$-colourings of the $n$-vertex path using systematic scan Markov
chains. An $H$-colouring of the $n$-vertex path is a \emph{homomorphism}
from the $n$-vertex path to some fixed graph $H$. We show that systematic
scan for $H$-colourings of the $n$-vertex path mixes in $O(\log n)$
scans for any fixed $H$. This is a significant improvement over the
previous bound on the mixing time which was $O(n^{5})$ scans. Furthermore
we show that for a slightly more restricted family of $H$ (where
any two vertices are connected by a 2-edge path) systematic scan also
mixes in $O(\log n)$ scans for \emph{any scan order}. Finally, for
completeness, we show that a random update Markov chain mixes in $O(n\log n)$
updates for any fixed $H$, improving the previous bound on the mixing
time from $O(n^{5})$ updates. 
\end{abstract}

\section{Introduction}

\begin{figure}

\caption{The graph describing the \emph{independent sets} model. Sites assigned
colour $a$ are ``out'' and sites assigned $b$ are ``in''.\label{fig:IS}}

\begin{centering}\includegraphics{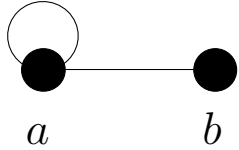}\par\end{centering}
\end{figure}

Many combinatorial problems are of interest to computer scientists
both in their own right and due to their natural applications to statistical
physics. Such problems can often be studied by considering \emph{homomorphisms}
from the graph of interest $G$ to some fixed graph $H$. This is
known as an $H$-colouring of $G$. The vertices of $H$ correspond
to colours and the edges of $H$ specify which colours are allowed
to be adjacent in an $H$-colouring of a graph. Let $H=(C,E)$ by
any fixed graph. We will refer to $C$ as the set of \emph{colours}
(in the literature it is often referred to as the set of \emph{spins}).
Formally an $H$-colouring of a graph $G=(V,E_{G})$ is a function
$h:V\rightarrow C$ such that $(h(v),h(u))\in E$ for all edges $(v,u)\in E_{G}$
of $G$. For example if $H$ is the graph in Figure~\ref{fig:IS}
then $C=\{ a,b\}$ and sites (in order to be consistent with existing
literature, e.g. Weitz~\cite{dror_combinatorial}, we will refer
to elements of $V$ as sites throughout this paper) assigned colour
$a$ in an $H$-colouring of $G$ are permitted to be adjacent to
sites assigned both $a$ and $b$, but sites assigned colour $b$
can only be adjacent to sites assigned colour $a$. 

\begin{figure}

\caption{The graph describing the \emph{Beach} model.\label{fig:Beach}}

\begin{centering}\includegraphics{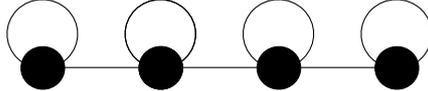}\par\end{centering}
\end{figure}

\begin{figure}

\caption{The graph describing the 4-particle \emph{Widom-Rowlinson} model.\label{fig:WR}}

\begin{centering}\includegraphics{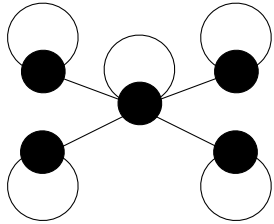}\par\end{centering}
\end{figure}

Due to the applicability of $H$-colourings to models in statistical
physics, and for ease of analysis, $H$-colouring problems are often
studied by restricting attention to a specific graph $H$. We now
give a few examples of special cases of $H$ that correspond to important
$H$-colouring problems. $H$-colourings using the graph $H$ from
Figure~\ref{fig:IS} correspond to \emph{independent set} configurations
of a graph where sites assigned colour $a$ are ``out'' and sites
assigned $b$ are ``in'' the independent set. It is usual to assign
weight $1$ to vertex $a$ and some positive weight $\lambda>0$ to
vertex $b$ in $H$. Independent sets (also known as the \emph{hard-core
lattice gas model} when using the weighted setting) is one of the
most commonly studied type of $H$-colourings in theoretical computer
science. Another well-studied case is when $H$ is the $q$-clique,
in which case $H$-colourings correspond to proper $q$-colourings
of the underlying graph. A proper $q$-colouring is a configuration
where no two adjacent sites are permitted to be assigned the same
colour. It is worth noting that proper $q$-colourings correspond
to the \emph{$q$}-state \emph{anti-ferromagnetic Potts model} at
zero temperature which is a well-studied model in statistical physics.
\emph{}Other well-known examples include the \emph{Beach} model introduced
by Burton and Steif~\cite{beach} and the $q$-particle Widom-Rowlinson
due to Widom and Rowlinson~\cite{WR}. The graph corresponding to
the Beach model is shown in Figure~\ref{fig:Beach}. The Beach model
was originally introduced as an example of a physical system, with
underlying graph $\Z^{d}$, which exhibits more than a single measure
of maximal entropy when $d>1$. The $q$-particle Widom-Rowlinson
model is a model of gas consisting of $q$ types of particles that
are not allowed to be adjacent to each other. The graph corresponding
to the $q=4$ case is shown in Figure~\ref{fig:WR} where the center
vertex represents empty sites and each remaining vertex represents
a particle.


The problem of determining whether a graph has an $H$-colouring for
a specific $H$ has been well-studied and Hell and Ne\v set\v ril~\cite{H-col-NP}
gave a complete characterisation of graphs $H$ for which this problem
is NP-complete. In particular, they showed that if $H$ has a loop
or is bipartite then the problem is in P, and that the problem is
NP-complete for any other fixed $H$. A complete dichotomy is also
known for the problem of counting the number of $H$-colourings. This
is due to Dyer and Greenhill~\cite{H-col-count} who showed that
if $H$ has at least one so-called \emph{nontrivial component} then
the counting problem is is \#P-complete. Otherwise it is in P. A trivial
component is a connected component which is either a complete graph
with all loops present, or a complete bipartite graph with no loops
present. They furthermore showed that the same dichotomy holds even
when the underlying graph is of bounded degree, which is significant
since in many physical applications the underlying graph tends to
be of low degree. Interestingly the above characterisation for the
decision problem does not hold for bounded degree graphs as shown
by Galluccio, Hell and Ne\v set\v ril~\cite{H-col-bounded}. Despite
the hardness of exactly counting the number of $H$-colourings it
remains possible to \emph{approximately} count the number of $H$-colourings
of a graph as long as it is possible to sample efficiently from the
(near) uniform distribution of $H$-colourings. This is due to a general
counting-to-sampling reduction of Dyer, Goldberg and Jerrum~\cite{H-col-sample}
that holds for any fixed $H$ and any underlying graph.

Sampling from the uniform distribution of $H$-colourings of a graph,
which for this discussion we will denote by $\pi$, is a challenging
task and some results about the complexity thereof are known. Goldberg,
Kelk and Paterson~\cite{kelk} have shown that, if $H$ has no nontrivial
components, then the sampling problem is intractable in a complexity-theoretic
sense. That is, they prove that there is unlikely to be any \emph{Polynomial
Almost Uniform Sampler} for $H$-colourings by reducing the problem
of sampling from the (near) uniform distribution of $H$-colourings
to the problem of counting independent sets in bipartite graphs, which
in turn is complete for a logically-defined subclass of \#P (see Dyer,
Goldberg, Greenhill and Jerrum~\cite{relative} for results about
this complexity class). This does, however, not rule out the possibility
of sampling from the uniform distribution of $H$-colourings of more
restricted graphs, such as the $n$-vertex path, as we will be focusing
on in this paper. This sampling task may be carried out by simulating
some suitable random dynamics converging to $\pi$. Ensuring that
a dynamics converges to $\pi$ is generally straightforward, but obtaining
good upper bounds on the number of steps required for the dynamics
to become sufficiently close to $\pi$ is a much more difficult problem.
Due to a lack of theoretical convergence results, scientists conducting
experiments by simulating such dynamics are at times forced to ``guess''
(using some heuristic methods) the number of steps required for their
dynamics to be sufficiently close to the desired distribution. See
for example Cowles and Carlin~\cite{MCMCDiagnostics} for a comprehensive
review of some diagnostic tools used to empirically determine these
convergence rates. By establishing rigorous bounds on the convergence
rates (\emph{mixing time}) of these dynamics computer scientists can
provide underpinnings for this type of experimental work and also
allow a more structured approach to be taken. 

Analysing the mixing time of Markov chains for $H$-colouring problems
is a well-studied area in theoretical computer science. There is a
substantial body of literature concerned with inventing Markov chains
for sampling from the uniform distribution of $H$-colourings of graphs
and providing bounds on their mixing times. When an $H$-colouring
corresponds to a proper $q$-colouring of graph with maximum vertex-degree
$\Delta$ then Jerrum~\cite{jerrum_simple}, and independently Salas
and Sokal~\cite{salas-sokal}, showed that a simple Markov chain
mixes in $O(n\log n)$ updates when $q>2\Delta$. This Markov chain
makes transitions by selecting a site $v$ and a colour $c$ uniformly
at random, and then recolouring site $v$ to $c$ if doing so results
in a proper $q$-colouring of the graph. By considering a more complicated
Markov chain Vigoda~\cite{vigoda} was able to weaken the restriction
on $q$ to $q>(11/6)\Delta$ colours being sufficient for proving
mixing in $O(n\log n)$ updates. This remains the least number of
colours required for mixing of a Markov chain on general graphs, however
the number of colours can be further reduced for special graphs. For
example, when the underlying graph is the square grid then Goldberg,
Martin and Paterson~\cite{ssm} gave a hand-proof that $q=7$ colours
are sufficient for mixing in $O(n\log n)$ updates by proving a condition
called ``strong spacial mixing''. Achlioptas, Molloy, Moore and van
Bussel~\cite{grid_ach} showed that $q=6$ colours are sufficient
for a Markov chain for proper colourings of the grid to mix in $O(n\log n)$
updates using a computational proof. As a final example for proper
$q$-colourings Martinelli, Sinclair and Weitz~\cite{tree} showed
that $q=\Delta+2$ colours are sufficient for $O(n\log n)$ mixing
when the underlying graph is a tree, improving a similar result by
Kenyon, Mossel and Peres~\cite{kenyon_tree}. When $H$ correspond
to independent set configurations of a graph with parameter $\lambda$
(that is, the vertex labeled $b$ in Figure~\ref{fig:IS} is assigned
some positive weight $\lambda$ and $a$ has weight 1) then $\lambda<\frac{2}{\Delta-2}$
is sufficient for $O(n\log n)$ mixing as shown by Dyer and Greenhill~\cite{ind_set}
and independently Luby and Vigoda~\cite{luby-vigoda-indset-improve}
(although the latter result is restricted to triangle-free graphs).
When $\Delta\leq4$ these results include the $\lambda=1$ case which
is of special interest to computer scientists since it corresponds
to sampling from the uniform distribution on independent sets of the
graph. Weitz~\cite{dror_ind_set} has recently improved the condition
on $\lambda$ to $\lambda<(\Delta-1)^{\Delta-1}/(\Delta-2)^{\Delta}$
which notably includes the $\lambda=1$ case for $\Delta=5$. When
$\Delta\geq6$ and $\lambda=1$ then Dyer, Frieze and Jerrum~\cite{ind_set_sparse}
have shown that there exists a bipartite graph $G_{0}$ such that
any so-called \emph{cautious} Markov chain on independent set configurations
of $G_{0}$ has (at least) exponential mixing time (in the number
of sites of $G_{0}$). A Markov chain is said to be cautious if it
is only allowed to change the state of a constant number of sites
at the time. This negative result was generalised to $H$-colourings
by Cooper, Dyer and Frieze~\cite{cooper}. Their result applies to
graphs $H$ that are either bipartite or have at least one loop present,
and is not a complete graph with all loops present (observe that for
such an $H$ the decision problem is in P and the counting problem
is in \#P as discussed above). In particular this result guarantees
the existence of a $\Delta$-regular graph $G_{0}$ (with $\Delta$
depending on $H$) such that any cautious Markov chain on the set
of $H$-colourings of $G_{0}$, and with uniform stationary distribution,
has a mixing time that is at least exponential in the number of sites
of $G_{0}$. 

While much is now understood about the mixing times of Markov chains,
the types of Markov chains frequently studied by computer scientists
fall under a family of Markov chains that we call \emph{random update
Markov chains}. We say that a Markov chain on the set of $H$-colourings
of a graph is a random update Markov chain when one step of the the
process consists of randomly selecting a set of sites (often a single
site) and updating the colours assigned to those sites according to
some well-defined distribution induced by $\pi$. The mixing time
of a random update Markov chain is measured in the number of \emph{updates}
required in order for the Markov chain to be sufficiently close to
$\pi$. We point out that all the positive results described above
are for random update Markov chains. An alternative to random update
Markov chains is to construct a Markov chain that cycles through and
updates the sites (or subsets of sites) in a deterministic order.
We call this a \emph{systematic scan Markov chain} (or \emph{systematic
scan} for short). \emph{}Systematic scan may be more intuitively appealing
in terms of implementation, however until recently this type of dynamics
has largely resisted analysis when applied to $H$-colouring problems.
Perhaps some of the first analyses of systematic scan were due to
Amit~\cite{amit} and Diaconis and Ram~\cite{diaconis-scan} who
respectively studied systematic scan in the context of sampling from
multivariate Gaussian distributions and generating random elements
of a finite group. The mixing time of a systematic scan Markov chain
is measured in the number of \emph{scans} of the graph required to
be sufficiently close to $\pi$ and throughout this paper it holds
that one scan takes $O(n)$ updates where $n$ is the number of sites
of the graph. It is important to note that systematic scan remains
a random process since the method used to update the colour assigned
to the selected set of sites is a randomised procedure drawing from
some well-defined distribution induced by $\pi$. This paper is concerned
with sampling from the uniform distribution of $H$-colourings of
the $n$-vertex path using systematic scan Markov chains. 

Only few results providing bounds on the mixing time of systematic
scan Markov chains for sampling from the uniform distribution of $H$-colourings
exist in the literature and almost all of them focus on proper $q$-colourings
of bounded degree graphs. For general graphs systematic scan is known
to mix in $O(\log n)$ scans whenever $q\geq2\Delta$, where $\Delta$
is the maximum vertex-degree of the graph, by updating both end-points
of an edge in each move. This is due to a recent result by Pedersen~\cite{dobrushin_ee}
which improves the polynomial in the $q=2\Delta$ case from a result
of Dyer, Goldberg and Jerrum~\cite{dobrushin_scan} that is obtained
by updating one site at the time. If the underlying graph is bipartite
then a systematic scan mixes in $O(\log n)$ scans whenever $q>f(\Delta)$
where $f(\Delta)\to\beta\Delta$ as $\Delta\to\infty$ and $\beta\approx1.76$.
This result is obtained by a careful construction of the metric used
in the coupling construction and is due to Bordewich, Dyer and Karpinski~\cite{metric}.
Furthermore, Dyer, Goldberg and Jerrum~\cite{systematic_scan} have
shown that a systematic scan for proper $3$-colourings of the $n$-vertex
path mixes in $\Theta(n^{2}\log n)$ scans when considering a systematic
scan which updates a single site at the time using the \emph{Metropolis}
update rule. In the same paper it is also shown that systematic scan
for $H$-colourings of the $n$-vertex path mixes in $O(n^{5})$ scans
for any fixed $H$ and that a random update Markov chain for $H$-colourings
of the $n$-vertex path mixes in $O(n^{5})$ updates. The authors
suggest, however, that both of these bounds are unlikely to be tight
and we will significantly improve them both in this paper.

In this paper we prove that systematic scan for $H$-colourings of
the $n$-vertex path mixes in $O(\log n)$ scans for any fixed graph
$H$, by updating a constant-size block of sites at each step. By
constant-size we mean that the number of sites contained in a block
is bounded independently of $n$. We do however allow the block-sizes
to depend on $H$ (since $H$ is a fixed graph). We will present two
different Markov chains in order to achieve this aim. In Section~\ref{sec:2path}
we show that if $H$ is a graph in which any pair of colours are connected
by a 2-edge path then a systematic scan mixes for any order of a set
of blocks, provided that the blocks are large enough. We will use
a recent result by Pedersen~\cite{dobrushin_ee}, which is based
on a technique known as \emph{Dobrushin uniqueness}, in order to establish
the mixing time of this Markov chain. In Section~\ref{sec:path_coupling}
we extend this result to all connected graphs $H$, although at the
expense of imposing a specific order on the scan. The proof of mixing
uses path coupling \cite{pathcoupling} in this case. Finally, for
completeness, we give a proof that a random update Markov chain for
$H$-colourings of the $n$-vertex path mixes in $O(n\log n)$ updates
for any fixed graph $H$. This result is presented in Section~\ref{sec:random-update}.

\subsection{Preliminaries and statement of results}

Consider a fixed (and connected) graph $H=(C,E)$ with maximum vertex-degree
$\Delta_{H}$. Let $C=\{1,\dots,q\}$ be referred to as the set of
colours. Also let $V=\{1,\ldots,n\}$ be the set of sites of the $n$-vertex
path and in particular let $V_{1}$ be the set of sites with odd indices
and $V_{2}$ the set of sites with even indices. We formally say that
an $H$-colouring of the $n$-vertex path is a function $h$ from
$V$ to $C$ such that $(h(i),h(i+1))\in E$ for all $i\in V\setminus\{ n\}$.
Let $\Omega^{+}$ be the set of all configurations (all possible assignments
of colours to the sites) of the $n$-vertex path and $\Omega$ be
the set of all $H$-colourings of the $n$-vertex path for the given
$H$. Define $\pi$ to be the uniform distribution on $\Omega$. If
$x\in\Omega^{+}$ is a configuration and $j\in V$ is a site then
$x_{j}$ denotes the colour assigned to $j$ in configuration $x$
and for any set $\Lambda\subseteq V$ let $x_{\Lambda}=\bigcup_{v\in\Lambda}\{ x_{v}\}$
be the set of colours assigned to sites in $\Lambda$. For colours
$c,d\in C$ and an integer $l$ let $D_{c,d}^{(l)}$ be the uniform
distribution on $H$-colourings of the region of consecutive sites
$L=\{ v_{1},\dots,v_{l}\}\subset V$ consistent with site $v_{1}$
being adjacent to a site $i\in V\setminus L$ assigned colour $c$
and site $v_{l}$ being adjacent to a site in $V\setminus L$ assigned
colour $d$. Also let $D_{c,d}^{(l)}(v_{j})$ be the distribution
on the colour assigned to site $v_{j}$ induced by $D_{c,d}^{(l)}$.
Observe that for $s<l$ \[
\left[D_{c,d}^{(l)}\mid v_{1}=c_{1},\dots,v_{s}=c_{s}\right]=D_{c_{s},d}^{(l-s)}\]
 where $D_{c,d}^{(l)}\mid v_{1}=c_{1},\dots,v_{s}=c_{s}$ is the uniform
distribution on $H$-colourings of $L$ conditioned on site $v_{1}$
being assigned colour $c_{1}$, $v_{2}$ colour $c_{2}$ and so on
until $v_{s}$ being assigned colour $c_{s}$. 

Let $\M$ be any ergodic Markov chain with state space $\Omega$ and
transition matrix $P$. By classical theory (see e.g. Aldous~\cite{aldous_walks})
$\M$ has a unique stationary distribution, which we will denote $\pi$.
The mixing time from an initial configuration $x\in\Omega$ is the
number of steps, that is applications of $P$, required for $\M$
to become sufficiently close to $\pi$. Formally the mixing time of
$\M$ from an initial configuration $x\in\Omega$ is defined, as a
function of the deviation $\epsilon$ from stationarity, by \[
\Mix_{x}(\M,\epsilon)=\min\{ t>0\;:\;\dtv(P^{t}(x,\cdot),\pi)\leq\epsilon\},\]
where \[
\dtv(\theta_{1},\theta_{2})=\frac{1}{2}\sum_{i}|\theta_{1}(i)-\theta_{2}(i)|=\max_{A\subseteq\Omega^{+}}|\theta_{1}(A)-\theta_{2}(A)|\]
is the total variation distance between two distributions $\theta_{1}$
and $\theta_{2}$ on $\Omega$. The mixing time $\Mix(M,\epsilon)$
of $\M$ is then obtained by maximising over all possible initial
configurations \[
\Mix(\M,\epsilon)=\max_{x\in\Omega}\Mix_{x}(\M,\epsilon).\]
We say that $\M$ is \emph{rapidly mixing} if the mixing time of $\M$
is polynomial in $n$ and $\log(\epsilon^{-1})$.

We study Markov chains that perform \emph{heat-bath} moves on a constant
number of sites at the time. For any configuration $x\in\Omega^{+}$
and subset of sites $\Lambda\subseteq V$ we let $\Omega_{\Lambda}(x)$
be the set of configurations where the colours assigned to the endpoints
of each edge containing a site in $\Lambda$ are also adjacent in
$H$. A heat-bath move on $\Lambda$ starting from configuration $x$
is performed by drawing a new configuration from the uniform distribution
on $\Omega_{\Lambda}(x)$. We would normally let $\Omega$ be the
state space of our Markov chains, however, if $H$ is bipartite then
we encounter a minor technical difficulty because the Markov chain
may not be ergodic. We overcome this ergodicity issue by partitioning
the state space as follows. If $C_{1}$ and $C_{2}$ are the colour
classes of $H$ then $\Omega_{1}=\{ x\in\Omega:x_{1}\in C_{1}\}$
is the set of $H$-colourings where the first site of the path is
assigned a colour from $C_{1}$. Observe that in fact each site in
$V_{1}$ is assigned a colour from $C_{1}$ and each site in $V_{2}$
is assigned a colour from $C_{2}$. Similarly $\Omega_{2}=\{ x\in\Omega:x_{1}\in C_{1}\}$
is the set of $H$-colourings where the first site is assigned a colour
from $C_{2}$. Intuitively, $\Omega_{1}$ and $\Omega_{2}$ are the
two connected components of $\Omega$ and we will show (Lemma~\ref{lem:ergodic-H-col-scan})
that the constructed Markov chains are ergodic on either $\Omega_{1}$
or $\Omega_{2}$. To see that $\Omega_{1}\cup\Omega_{2}$ contain
all $H$-colourings of the $n$-vertex path it is enough to observe
that if $x\in\Omega$ then any pair of adjacent sites of the $n$-vertex
path must be assigned colours from opposite colour classes of $H$
in $x$. We let $\state$ be the relevant state space of the Markov
chains in order to ensure ergodicity. In particular, if $H$ is non-bipartite
then $\state=\Omega$. Otherwise $H$ is bipartite and we let $\state$
be one of $\Omega_{1}$ and $\Omega_{2}$. This is the same partition
used by Dyer et al. in \cite{systematic_scan}. See also Cooper et
al.~\cite{cooper} for a discussion of this issue in the context
of $H$-colourings.

We are now ready to formally define the systematic scan Markov chains
we will study in this paper and state our theorems. Let $l_{1}=\lceil\Delta_{H}^{2}\log(\Delta_{H}^{2}+1)\rceil+1$.
Then let $\{\Theta_{1},\dots,\Theta_{m_{1}}\}$ be any set of $m_{1}=\lceil n/l_{1}\rceil$
blocks where each block consists of $l_{1}$ consecutive sites and
$\bigcup_{k=1}^{m_{1}}\Theta_{k}=V$. For each block $\Theta_{k}$
we define $\Pk$ to be the transition matrix on the state space $\state$
for performing a heat-bath move on $\Theta_{k}$. 

\begin{defn}
For any integer $n$ we let $\HscanDob$ be the systematic scan Markov
chain, on the state space $\state$, with transition matrix $\Pi_{k=1}^{m_{1}}\Pk$. 
\end{defn}
It is worth pointing out that the following result holds for \emph{any}
\emph{order} of the blocks, as is the case for all results obtained
by Dobrushin uniqueness (see e.g. Dyer et al~\cite{dobrushin_scan}).
In Section \ref{sec:2path} we will use a recent result by Pedersen~\cite{dobrushin_ee}
to prove the following theorem.

\begin{thm}
\label{thm:Hcol_2path} Let $H$ be a fixed connected graph and consider
the systematic scan Markov chain $\HscanDob$ on the state space $\state$.
Suppose that $H$ is a graph in which every two sites are connected
by a $2$-edge path. Then mixing time of $\HscanDob$ is \[
\Mix(\HscanDob,\epsilon)\leq\Delta_{H}^{2}(\Delta_{H}^{2}+1)\log(n\epsilon^{-1})\]
scans of the $n$-vertex path. This corresponds to $O(n\log n)$ updates
by the construction of the set of blocks.
\end{thm}
\begin{rem*}
Note that that each $H$ for which Theorem~\ref{thm:Hcol_2path}
is valid is non-bipartite so $\state=\Omega$.
\end{rem*}

\begin{rem*}
Several well known graphs satisfy the condition of Theorem \ref{thm:Hcol_2path},
for example Widom-Rowlinson configurations, independent set configurations
and proper $q$-colourings for $q\geq3$. The fact that an $H$ corresponding
to 3-colourings satisfies the condition of the theorem is particularly
interesting since a lower bound of $\Omega(n^{2}\log n)$ scans for
single site systematic scan on the path is proved in Dyer at al.~\cite{systematic_scan}.
This means that using a simple single site coupling cannot be sufficient
to establishing Theorem~\ref{thm:Hcol_2path} for any family of $H$
including 3-colourings and hence we have to use block updates. 
\end{rem*}
While many natural $H$-colouring problems belong to the family covered
by Theorem~\ref{thm:Hcol_2path}, others (e.g. Beach configurations)
are not included. We go on to show that systematic scan mixes in $O(\log n)$
scans for any fixed graph $H$ by placing more strict restrictions
on the construction of the blocks and the order of the scan. Let $s=4q+1$,
$\beta=\lceil\log(2sq^{s}+1)\rceil q^{s}$ and $l_{2}=2\beta s$.
For any integer $n$ consider the following set of $m_{2}+1=\lfloor2n/l_{2}\rfloor$
blocks $\{\Theta_{0},\dots,\Theta_{m_{2}}\}$ where \[
\Theta_{k}=\{ k\beta s+1,\dots,\min((k+2)\beta s,n)\}.\]
We observe that $\bigcup_{k=0}^{m_{2}}\Theta_{k}=V$ by construction
of the set of blocks. Furthermore note that the size of $\Theta_{m_{2}}$
is at least $\beta s$ and that the size of every other block is $l_{2}$. 

\begin{defn}
For any integer $n$ we let $\HscanPC$ be the systematic scan Markov
chain, on the state space $\state$, which performs a heat-bath move
on each block in the order $\Theta_{0},\dots,\Theta_{m_{2}}$. 
\end{defn}
In Section \ref{sec:path_coupling} we will use path coupling \cite{pathcoupling}
to prove the following theorem, which improves the mixing time from
the corresponding result in Dyer et al.~\cite{systematic_scan} from
$O(n^{5})$ scans to $O(\log n)$ scans.

\begin{thm}
\label{thm:Hcol_line_path} Let $H$ be any fixed connected graph
and consider the systematic scan Markov chain $\HscanPC$ on the state
space $\state$. The mixing time of $\HscanPC$ is \[
\Mix(\HscanPC,\epsilon)\leq(4sq^{s}+2)\log(n\epsilon^{-1})\]
scans of the $n$-vertex path. This corresponds to $O(n\log n)$ updates
by the construction of the set of blocks.
\end{thm}
\begin{rem*}
It is worth remarking at this point that Theorem~\ref{thm:Hcol_line_path}
eclipses Theorem~\ref{thm:Hcol_2path} in the sense that it shows
the existence of a systematic scan for a broader family of $H$ than
Theorem~\ref{thm:Hcol_2path} but with the same (asymptotic) mixing
time. The result stated as Theorem~\ref{thm:Hcol_2path} however
remains interesting in its own right since it applies to any order
of the scan. Following the proof of Theorem~\ref{thm:Hcol_2path}
we will discuss (Observation~\ref{obs:noextension}) the obstacles
one encounters when attempting to extend Theorem~\ref{thm:Hcol_2path}
to a larger family of $H$ using the same method of proof.
\end{rem*}
For completeness we finally consider a random update Markov chain
for $H$-colourings of the $n$-vertex path. Let $\gamma=2q^{s}+1$
and define the following set of $n+s\gamma-1$ blocks, which is constructed
such that each site is contained in exactly $s\gamma$ blocks \[
\Theta_{k}=\begin{cases}
\{ k,\dots,\min(k+s\gamma-1,n)\} & \textnormal{when }k\in\{1,\dots,n\}\\
\{1,\dots,n+s\gamma-k\} & \textnormal{when }k\in\{ n+1,\dots,n+s\gamma-1\}.\end{cases}\]

\begin{defn}
For any integer $n$ we let $\Hrandom$ be the random update Markov
chain, on the state space $\state$, which at each step selects a
block uniformly at random and performs a heat-bath move on it. 
\end{defn}
In Section \ref{sec:random-update} we will use path coupling \cite{pathcoupling}
to prove the following theorem, which improves the mixing time from
the corresponding result in Dyer et al.~\cite{systematic_scan} from
$O(n^{5})$ updates to $O(n\log n)$ updates.

\begin{thm}
\label{thm:random_Hcolor}Let $H$ be any fixed connected graph and
consider the random update Markov chain $\Hrandom$ on the state space
$\state$. The mixing time of $\Hrandom$ is \[
\Mix(\Hrandom,\epsilon)\leq\frac{(n+2sq^{s}+s-1)\log(n\epsilon^{-1})}{s}\]
block-updates. This corresponds to $O(n\log n)$ updates since the
size of each block is at most $s\gamma=O(1)$.
\end{thm}

\subsection{Review of proof techniques}

We now briefly introduce the techniques we will use to bound the mixing
time of the above Markov chains. For technical reasons we extend the
state space of the Markov chains as follows. Let $\Omega_{1}^{+}$
be the set of configurations where each site in $V_{1}$ is assigned
a colour from $C_{1}$ and each site in $V_{2}$ is assigned a colour
from $C_{2}$ (recall that $C_{1}$ and $C_{2}$ are the colour classes
of $H$). Similarly, $\Omega_{2}^{+}$ is the set of configurations
where each site in $V_{1}$ is assigned a colour from $C_{2}$ and
each site in $V_{2}$ is assigned a colour from $C_{1}$. Formally
\[
\Omega_{1}^{+}=\{ x\in\Omega^{+}:x_{V_{1}}\subseteq C_{1},x_{V_{2}}\subseteq C_{2}\}\]
and\[
\Omega_{2}^{+}=\{ x\in\Omega^{+}:x_{V_{1}}\subseteq C_{2},x_{V_{2}}\subseteq C_{1}\}.\]
We then extend the state space of the Markov chains to $\state^{+}$
where $\state^{+}=\Omega^{+}$ if $H$ is not bipartite and $\state^{+}$
is one of $\Omega_{1}^{+}$ or $\Omega_{2}^{+}$ when $H$ is bipartite.
The extended Markov chains make the same transitions as the original
Markov chains on configurations in $\state$ and hence the extended
chains do not make transitions from configurations in $\state$ to
configurations outside $\state$. The stationary distributions of
the extended chains are uniform over the configurations in $\state$
and zero elsewhere. This approach is standard and the mixing times
of the original chains are bounded above by the mixing time of corresponding
chain on the extended state space. 

For each site $j\in V$, let $S_{j}$ denote the set of pairs $(x,y)\in\state^{+}\times\state^{+}$
of configurations that only differ on the colour assigned to site
$j$, that is $x_{i}=y_{i}$ for all $i\neq j$. Also let $S=\bigcup_{j\in V}S_{j}$
be the set of all such pairs of configurations.

\subsubsection{Dobrushin uniqueness}

We will make use of a recent result by Pedersen~\cite{dobrushin_ee}
to prove Theorem \ref{thm:Hcol_2path} by bounding the influence \emph{on}
a site. For completeness we now summarise this result and at the same
time point out how the construction of $\HscanDob$ ensures that all
required properties are satisfied. First note from the remark after
Theorem~\ref{thm:Hcol_2path} that each $H$ that we consider is
not bipartite and so $\state^{+}=\Omega^{+}$. Suppose that $\{\Theta_{1},\dots,\Theta_{m}\}\subseteq V$
is a set of $m$ blocks such that $\bigcup_{k=1}^{m}\Theta_{k}=V$
and that each block $\Theta_{k}$ is associated with a transition
matrix $\Pk$ on the state space $\Omega^{+}$. For any configuration
$x\in\Omega^{+}$, $\Pk(x,\cdot)$ denotes the distribution on configurations
obtained from applying $\Pk$ to $x$. Recall from the definition
of $\HscanDob$ that the set of blocks covers $V$ as required and
that each transition matrix $\Pk$ represents performing a heat-bath
move on $\Theta_{k}$. It is furthermore required that each transition
matrix $\Pk$ satisfies the following two properties.

\begin{enumerate}
\item If $\Pk(x,y)>0$ then $x_{i}=y_{i}$ for each $i\in V\setminus\Theta_{k}$,
and
\item the distribution $\pi$ on $\Omega^{+}$ is invariant with respect
to $\Pk.$
\end{enumerate}
Pedersen \cite{dobrushin_ee} points out that if $\Pk$ is the transition
matrix performing a heat-bath move on $\Theta_{k}$ and $\pi$ is
the uniform distribution on $\Omega$, as they both are in the case
of $\HscanDob$, then both of these properties are satisfied. These
properties ensure that the stationary distribution of any systematic
scan Markov chain with transition matrix $\Pi_{k=1}^{m}\Pk$ is $\pi$. 

We are now ready to define the parameter $\alpha$ denoting the influence
on a site. For any pair of configurations $(x,y)\in S_{i}$ let $\Psi_{k}(x,y)$
be a coupling of the distributions $\Pk(x,\cdot)$ and $\Pk(y,\cdot)$.
We remind the reader that a coupling of two distributions $\pi_{1}$
and $\pi_{2}$ on state space $\Omega^{+}$ is a joint distribution
on $\Omega^{+}\times\Omega^{+}$ such that the marginal distributions
are $\pi_{1}$ and $\pi_{2}$. We let $(x^{\prime},y^{\prime})\in\Psi_{k}(x,y)$
denote that the pair of configurations $(x^{\prime},y^{\prime})$
is drawn from the coupling $\Psi_{k}(x,y)$. We then let \[
\rho_{i,j}^{k}=\max_{(x,y)\in S_{i}}\Prob_{(x^{\prime},y^{\prime})\in\Psi_{k}(x,y)}(x_{j}^{\prime}\neq y_{j}^{\prime})\]
be the influence of site $i$ on site $j$ under $\Theta_{k}$. The
influence of $i$ on $j$ is thus the probability that site $j\in\Theta_{k}$
is assigned a different colour in a pair of configurations drawn from
the coupling $\Psi_{k}(x,y)$ where $x$ and $y$ differ only on the
colour of site $i$. Finally the parameter $\alpha$ denoting the
influence on any site is defined as \[
\alpha=\max_{k}\max_{j\in\Theta_{k}}\sum_{i}\rho_{i,j}^{k}.\]

\begin{rem*}
Pedersen~\cite{dobrushin_ee} actually defines $\alpha$ with a positive
weight assigned to each site of the graph, however as we will not
use the weights in our proof they are omitted from the above definition. 
\end{rem*}
The following theorem bounds the mixing time of a systematic scan
Markov chain $\scan$ with transition matrix $\Pi_{k=1}^{m}\Pk$.
It is worth pointing out that, since the proof makes use of Dobrushin
uniqueness, this upper-bound on the mixing time holds for any order
of the blocks.

\begin{thm}
[Theorem 2, Pedersen \cite{dobrushin_ee}]\label{thm:main_d} If $\alpha<1$
then the mixing time of $\scan$ satisfies \[
\Mix(\scan,\epsilon)\leq\frac{\log(n\epsilon^{-1})}{1-\alpha}.\]

\end{thm}

\subsubsection{Path coupling}

In order to prove Theorems~\ref{thm:Hcol_line_path} and~\ref{thm:random_Hcolor}
we will make use of path coupling~\cite{pathcoupling} which is a
well-known, and by now standard, technique for proving rapid mixing
of Markov chains. The key idea of path coupling is to define a coupling
for pairs of \emph{adjacent} configurations where the set of all adjacent
configurations connects the state space. We will say that a pair of
configurations $x,y\in\state^{+}$ are adjacent if $(x,y)\in S$.
The path coupling machinery then extends the coupling to all pairs
of configurations in the state space. For completeness we show that
$S$ connects the state space $\state^{+}$.

\begin{lem}
The transitive closure of $S$ is the whole of $\state^{+}\times\state^{+}$.
\end{lem}
\begin{proof}
Recall that $S=\bigcup_{i\in V}S_{j}$ where $S_{j}\subseteq\state^{+}\times\state^{+}$
is the set of pairs $(x,y)\in\state^{+}\times\state^{+}$ of configurations
that differ only on the colour assigned to site $j$. To establish
the lemma it is sufficient to, for any pair of configurations $(x,y)\in\state^{+}\times\state^{+}$,
to construct a path $x=z^{0},z^{1},\dots,z^{n}=y$ such that $(z^{j-1},z^{j})\in S_{j}$
for each $j\in\{1,\dots,n\}$. We define $z^{j}$ for $j\in\{1,\dots,n\}$
as follows \[
z_{i}^{j}=\begin{cases}
y_{i} & \textnormal{for }1\leq i\leq j\\
x_{i} & \textnormal{for }j<i\leq n.\end{cases}\]
Informally, configuration $z^{j}$ agrees with configuration $y$
from site $1$ to $j$ and with configuration $x$ from site $j+1$
to $n$. 

By definition of the configurations $z^{0},\dots,n^{n}$ it follows
that $z^{j-1}$ and $z^{j}$ only differ on the colour assigned to
site $j$ for each $j\in\{1,\dots,n\}$. Hence we only need to check
that $z^{j}\in\state^{+}$ for each $j$. If $H$ is non-bipartite
then $\state^{+}=\Omega^{+}$ so $z^{j}\in\state^{+}$ for each $j\in\{1,\dots,n\}$.
If $H$ is bipartite then $\state^{+}$ is one of $\Omega_{1}^{+}$
or $\Omega_{2}^{+}$. Suppose without loss of generality that $\state^{+}=\Omega_{1}^{+}$.
Then for each $j\in\{1,\dots n\}$ it holds by definition of $\Omega_{1}^{+}$
that the colours $x_{j}$ and $y_{j}$ must be from the same colour
class of $H$ and hence have $z^{j}\in\Omega_{1}^{+}$.
\end{proof}
Finally note that $\Ham(x,y)=1$ for any $(x,y)\in S$ where $\Ham(x,y)$
denotes the Hamming distance between configurations $x$ and $y$.
The following theorem is sufficient for our needs in this paper, and
it is a special case of the general path coupling theorem proved by
Bubley and Dyer~\cite{pathcoupling}. 

\begin{thm}
[Bubley, Dyer \cite{pathcoupling}]\label{thm:path-coupling}For all
pairs $(x,y)\in S$ define a coupling $(x,y)\mapsto(x^{\prime},y^{\prime})$
of a Markov chain $\M$ on the state space $\state^{+}$. Suppose
that there exists a constant $0<\gamma<1$ such that $\E{\Ham(x^{\prime},y^{\prime})}\leq(1-\gamma)$
for all pairs $(x,y)\in S$. Then the mixing time of $\M$ satisfies\[
\Mix(\M,\epsilon)\leq\frac{\log(n\epsilon^{-1})}{\gamma}.\]

\end{thm}

\section{$H$-colourings on the path for a restricted family of $H$\label{sec:2path}}

Recall that $\Delta_{H}$ denotes the maximum vertex-degree of some
fixed graph $H$ and that $l_{1}=\lceil\Delta_{H}^{2}\log(\Delta_{H}^{2}+1)\rceil+1$.
The systematic scan Markov chain $\HscanDob$ on $\state$ has transition
matrix $\Pi_{k=1}^{m_{1}}\Pk$ where $\Pk$ is the transition matrix
for performing a heat-bath move on block $\Theta_{k}$ from a set
of $m_{1}=\lceil n/l_{1}\rceil$ size $l_{1}$ blocks covering the
$n$-vertex path. We will prove Theorem~\ref{thm:Hcol_2path}, namely
that $\HscanDob$ mixes in $O(\log n)$ scans when $H$ is a graph
in which any two colours are connected via a 2-edge path. We will
bound the mixing time of $\HscanDob$ by bounding the influence on
a site and begin by establishing some lemmas required to construct
the coupling needed in the proof of Theorem~\ref{thm:Hcol_2path}.

\begin{lem}
\label{lemma:first-coupling} Suppose that for any $c_{1},c_{2}\in C$
there is a 2-edge path in $H$ from $c_{1}$ to $c_{2}$. Then for
any $c_{1},c_{2},d\in C$ and integer $s'\geq2$ there exists a coupling
$\psi(D_{c_{1},d}^{(s')},D_{c_{2},d}^{(s')})$ of $D_{c_{1},d}^{(s')}$
and $D_{c_{2},d}^{(s')}$ such that 

\begin{enumerate}[{(i)}]
\item $\Prob_{(x^{\prime},y^{\prime})\in\psi(D_{c_{1},d}^{(s')},D_{c_{2},d}^{(s')})}(x_{v_{1}}^{\prime}\neq y_{v_{1}}^{\prime})\leq1-\frac{1}{\Delta_{H}^{2}}$
and 

\item $\Prob_{(x^{\prime},y^{\prime})\in\psi(D_{c_{1},d}^{(2)},D_{c_{2},d}^{(2)})}(x_{v_{2}}^{\prime}\neq y_{v_{2}}^{\prime})\leq1-\frac{1}{\Delta_{H}^{2}}.$
\end{enumerate}
\end{lem}
\begin{proof}
By the condition of the lemma there exists some $c^{\prime}\in C$
adjacent to both $c_{1}$ and $c_{2}$ in $H$. We prove the statement
by considering two cases on $s$.

First suppose that $s'=2$. By the condition of the lemma there is
some colour $d^{\prime}$ adjacent to both $c^{\prime}$ and $d$
in $H$. There are at most $\Delta_{H}^{2}$ valid $H$-colourings
of the sites $v_{1},v_{2}$ in either of the distributions $D_{c_{1},d}^{(2)}$
and $D_{c_{2},d}^{(2)}$, and hence the colouring $h$, which assigns
$c^{\prime}$ to $v_{1}$ and $d^{\prime}$ to $v_{2}$, has weight
at least $1/\Delta_{H}^{2}$ in both. We construct a coupling $\psi(D_{c_{1},d}^{(2)},D_{c_{2},d}^{(2)})$
such that \[
\Prob_{(x^{\prime},y^{\prime})\in\psi(D_{c_{1},d}^{(2)},D_{c_{2},d}^{(2)})}(x^{\prime}=y^{\prime}=h)\geq\frac{1}{\Delta_{H}^{2}}.\]
 The rest of the coupling is arbitrary. This gives the following bounds
on the disagreement probabilities at $v_{1}$ and $v_{2}$ \[
\Prob_{(x^{\prime},y^{\prime})\in\psi(D_{c_{1},d}^{(2)},D_{c_{2},d}^{(2)})}(x_{v_{1}}^{\prime}=y_{v_{1}}^{\prime})\geq\Prob_{(x^{\prime},y^{\prime})\in\psi(D_{c_{1},d}^{(2)},D_{c_{2},d}^{(2)})}(x_{v_{1}}^{\prime}=y_{v_{1}}^{\prime}=c^{\prime})\geq\frac{1}{\Delta_{H}^{2}}\]
 which establishes $(i)$ for $s'=2$ and \[
\Prob_{(x^{\prime},y^{\prime})\in\psi(D_{c_{1},d}^{(2)},D_{c_{2},d}^{(2)})}(x_{v_{2}}^{\prime}=y_{v_{2}}^{\prime})\geq\Prob_{(x^{\prime},y^{\prime})\in\psi(D_{c_{1},d}^{(2)},D_{c_{2},d}^{(2)})}(x_{v_{2}}^{\prime}=y_{v_{2}}^{\prime}=d^{\prime})\geq\frac{1}{\Delta_{H}^{2}}\]
 which establishes $(ii)$.

Now suppose $s'>2$. Let $adj(c)$ denote the set of colours adjacent
to $c$ in $H$ and $n_{k}$ the number of $H$-colourings on the
sites $v_{4},\dots,v_{s'}$ consistent with $v_{3}$ being assigned
colour $k\in C$ and $v_{s'}$ being adjacent to a site (out side
the block) coloured $d$. Also let $p_{c,k}$ be the number of $H$-colourings
of $v_{1},v_{2},v_{3}$ assigning colour $c$ to $v_{1}$ and $k$
to $v_{3}$ without regard to other sites. Finally let $z_{i}$ be
the number of $H$-colourings with positive measure in $D_{c_{i},d}^{(s')}$
and assume without loss of generality that $z_{1}\geq z_{2}$.

There are at most $\Delta_{H}$ colours available for each site in
the block which gives $p_{c,k}\leq\Delta_{H}$ for any $c,k\in C$
and hence \[
z_{1}=\sum_{c\in adj(c_{1})}\sum_{k\in C}p_{c,k}n_{k}\leq\Delta_{H}\sum_{c\in adj(c_{1})}\sum_{k\in C}n_{k}\leq\Delta_{H}^{2}\sum_{k\in C}n_{k}.\]
Now let $H(c^{\prime})$ be the set of all $H$-colourings with positive
measure in $D_{c_{1},d}^{(s')}$ that assign colour $c^{\prime}$
to site $v_{1}$. Let $h(c^{\prime})$ denote the size of this set.
Now $p_{c,k}\geq1$ for any $c,k\in C$ since there is a 2-edge path
in $H$ between any two colours and hence \[
h(c^{\prime})=\sum_{k\in C}p_{c^{\prime},k}n_{k}\geq\sum_{k\in C}n_{k}.\]

Observe that, for any $h\in H(c^{\prime})$, $h$ is at least as likely
in $D_{c_{2},d}^{(s')}$ as in $D_{c_{1},d}^{(s')}$ since we have
assumed $z_{1}\geq z_{2}$ without loss of generality. We construct
a coupling $\psi(D_{c_{1},d}^{(s')},D_{c_{2},d}^{(s')})$ of $D_{c_{1},d}^{(s')}$
and $D_{c_{2},d}^{(s')}$ in which for each $h\in H(c^{\prime})$
\[
\Prob_{(x^{\prime},y^{\prime})\in\psi(D_{c_{1},d}^{(s')},D_{c_{2},d}^{(s')})}(x^{\prime}=y^{\prime}=h)\geq\frac{1}{z_{1}}.\]
 The rest of the coupling is arbitrary. Hence \begin{align*}
\Prob_{(x^{\prime},y^{\prime})\in\psi(D_{c_{1},d}^{(s')},D_{c_{2},d}^{(s')})}(x_{v_{1}}^{\prime}=y_{v_{1}}^{\prime}) & \geq\sum_{h\in H(c^{\prime})}\Prob_{(x^{\prime},y^{\prime})\in\psi(D_{c_{1},d}^{(s')},D_{c_{2},d}^{(s')})}(x^{\prime}=y^{\prime}=h)\\
 & \geq\frac{h(c^{\prime})}{z_{1}}\\
 & \geq\frac{1}{\Delta_{H}^{2}}\end{align*}
using the bounds on $z_{1}$ and $h(c^{\prime})$. This completes
the proof.
\end{proof}
We then use Lemma~\ref{lemma:first-coupling} to bound the disagreement
probabilities at each site of of the block when a pair of configurations
are drawn from a recursively constructed coupling. 

\begin{lem}
\label{lemma:line_coup}Suppose that for any $c_{1},c_{2}\in C$ there
is a 2-edge path in $H$ from $c_{1}$ to $c_{2}$. Then for all $c_{1},c_{2},d\in C$
and integers $l^{\prime}\geq2$ there exists a coupling $\Psi(D_{c_{1},d}^{(l')},D_{c_{2},d}^{(l')})$
of $D_{c_{1},d}^{(l')}$ and $D_{c_{2},d}^{(l')}$ in which for $j\in\{1,\dots,l'-1\}$
\[
\Prob_{(x^{\prime},y^{\prime})\in\Psi(D_{c_{1},d}^{(l')},D_{c_{2},d}^{(l')})}(x_{v_{j}}^{\prime}\neq y_{v_{j}}^{\prime})\leq\left(1-\frac{1}{\Delta_{H}^{2}}\right)^{j}\]
 and \[
\Prob_{(x^{\prime},y^{\prime})\in\Psi(D_{c_{1},d}^{(l')},D_{c_{2},d}^{(l')})}(x_{v_{l}}^{\prime}\neq y_{v_{l}}^{\prime})\leq\left(1-\frac{1}{\Delta_{H}^{2}}\right)^{l'-1}.\]

\end{lem}
\begin{proof}
We recursively construct a coupling $\Psi(D_{c_{1},d}^{(l')},D_{c_{2},d}^{(l')})$
of $D_{c_{1},d}^{(l')}$ and $D_{c_{2},d}^{(l')}$ using the method
set out in Goldberg et al.~\cite{ssm} as follows. Firstly $l'=2$
is the base case and we use the coupling from Lemma \ref{lemma:first-coupling}.
For $l'\geq3$ we construct a coupling using the following two step
process. 
\begin{enumerate}
\item Couple $D_{c_{1},d}^{(l')}(v_{1})$ and $D_{c_{2},d}^{(l')}(v_{1})$
greedily to maximise the probability of assigning the same colour
to site $v_{1}$ in both distributions. 
\item If the same colour $c$ was chosen for $v_{1}$ in both distributions
in step 1 then the set of valid $H$-colourings of the remaining sites
are the same in both distributions. Hence the conditional distributions
$D_{c_{1},d}^{(l')}\mid v_{1}=c$ and $D_{c_{2},d}^{(l')}\mid v_{1}=c$
are the same and the rest of the coupling is trivial. Otherwise, for
all pairs $(c_{1}^{\prime},c_{2}^{\prime})$ of distinct colours recursively
couple $D_{c_{1},d}^{(l')}\mid v_{1}=c_{1}^{\prime}=D_{c_{1}',d}^{(l'-1)}$
and $D_{c_{2},d}^{(l')}\mid v_{1}=c_{2}^{\prime}=D_{c_{2}',d}^{(l'-1)}$
which is a sub problem of size $l'-1$. 
\end{enumerate}
This completes the coupling construction.

Now for $j\in\{1,\dots,l'-1\}$ we prove by induction that \begin{equation}
\Prob_{(x^{\prime},y^{\prime})\in\Psi(D_{c_{1},d}^{(l')},D_{c_{2},d}^{(l')})}(x_{v_{j}}^{\prime}\neq y_{v_{j}}^{\prime})\leq\left(1-\frac{1}{\Delta_{H}^{2}}\right)^{j}.\label{eq:2path-induct}\end{equation}
 The base case, $j=1$, follows from Lemma \ref{lemma:first-coupling}
since we couple the colour at site $v_{1}$ greedily to maximise the
probability of agreement at $v_{1}$ in the first step of the recursive
coupling. Now suppose that \eqref{eq:2path-induct} is true for $j-1$
then \begin{align*}
\Prob & _{(x^{\prime},y^{\prime})\in\Psi(D_{c_{1},d}^{(l')},D_{c_{2},d}^{(l')})}(x_{v_{j}}^{\prime}\neq y_{v_{j}}^{\prime})\\
 & =\sum_{c_{1}^{\prime},c_{2}^{\prime}}\Prob_{(x^{\prime},y^{\prime})\in\Psi(D_{c_{1},d}^{(l')},D_{c_{2},d}^{(l')})}(x_{v_{j-1}}^{\prime}=c_{1}^{\prime},y_{v_{j-1}}^{\prime}=c_{2}^{\prime})\\
 & \quad\times\Prob_{(x^{\prime},y^{\prime})\in\Psi(D_{c_{1},d}^{(l')}\mid v_{j-1}=c_{1}^{\prime},D_{c_{2},d}^{(l')}\mid v_{j-1}=c_{2}^{\prime})}(x_{v_{j}}^{\prime}\neq y_{v_{j}}^{\prime})\\
 & =\sum_{c_{1}^{\prime},c_{2}^{\prime}}\Prob_{(x^{\prime},y^{\prime})\in\Psi(D_{c_{1},d}^{(l')},D_{c_{2},d}^{(l')})}(x_{v_{j-1}}^{\prime}=c_{1}^{\prime},y_{v_{j-1}}^{\prime}=c_{2}^{\prime})\\
 & \quad\times\Prob_{(x^{\prime},y^{\prime})\in\Psi(D_{c_{1}^{\prime},d}^{(l'-j+1)},D_{c_{2}^{\prime},d}^{(l'-j+1)})}(x_{v_{1}}^{\prime}\neq y_{v_{1}}^{\prime})\\
 & \leq\sum_{c_{1}^{\prime},c_{2}^{\prime}}\Prob_{(x^{\prime},y^{\prime})\in\Psi(D_{c_{1},d}^{(l')},D_{c_{2},d}^{(l')})}(x_{v_{j-1}}^{\prime}=c_{1}^{\prime},y_{v_{j-1}}^{\prime}=c_{2}^{\prime})\left(1-\frac{1}{\Delta_{H}^{2}}\right)\\
 & \leq\left(1-\frac{1}{\Delta_{H}^{2}}\right)^{j}\end{align*}
 where the first inequality uses Lemma \ref{lemma:first-coupling}
and the second is the inductive hypothesis.

The $j=l'$ case is similar. \begin{align*}
\Prob & _{(x^{\prime},y^{\prime})\in\Psi(D_{c_{1},d}^{(l')},D_{c_{2},d}^{(l')})}(x_{v_{l}}^{\prime}\neq y_{v_{l}}^{\prime})\\
 & =\sum_{c_{1}^{\prime},c_{2}^{\prime}}\Prob_{(x^{\prime},y^{\prime})\in\Psi(D_{c_{1},d}^{(l')},D_{c_{2},d}^{(l')})}(x_{v_{l'-2}}^{\prime}=c_{1}^{\prime},y_{v_{l'-2}}^{\prime}=c_{2}^{\prime})\\
 & \quad\times\Prob_{(x^{\prime},y^{\prime})\in\Psi(D_{c_{1},d}^{(l')}\mid v_{l'-2}=c_{1}^{\prime},D_{c_{2},d}^{(l')}\mid v_{l'-2}=c_{2}^{\prime})}(x_{v_{l'}}^{\prime}\neq y_{v_{l'}}^{\prime})\\
 & =\sum_{c_{1}^{\prime},c_{2}^{\prime}}\Prob_{(x^{\prime},y^{\prime})\in\Psi(D_{c_{1},d}^{(l')},D_{c_{2},d}^{(l')})}(x_{v_{l'-2}}^{\prime}=c_{1}^{\prime}\land y_{v_{l'-2}}^{\prime}=c_{2}^{\prime})\\
 & \quad\times\Prob_{(x^{\prime},y^{\prime})\in\Psi(D_{c_{1}^{\prime},d}^{(2)},D_{c_{2}^{\prime},d}^{(2)})}(x_{v_{2}}^{\prime}\neq y_{v_{2}}^{\prime})\\
 & \leq\left(1-\frac{1}{\Delta_{H}^{2}}\right)^{l'-2}\left(1-\frac{1}{\Delta_{H}^{2}}\right)=\left(1-\frac{1}{\Delta_{H}^{2}}\right)^{l'-1}\end{align*}
 where the inequality uses Lemma \ref{lemma:first-coupling} and \eqref{eq:2path-induct}.
\end{proof}
We can then use the coupling constructed in Lemma~\ref{lemma:line_coup}
to construct a coupling $\Psi_{k}(x,y)$ of the distributions $\Pk(x,\cdot)$
and $\Pk(y,\cdot)$ for each pair of configurations $(x,y)\in S_{i}$.
We summarise the disagreement probabilities in this coupling in the
following corollary (of Lemma~\ref{lemma:line_coup}). 

\begin{cor}
\label{cor:Hcol_rho} For any sites $i,j\in V$ let $d(i,j)$ denote
the edge distance between them and suppose that for any $c,d\in C$
there exists a $2$-edge path in $H$ from $c$ to $d$. Then \begin{eqnarray*}
\rho_{i,j}^{k} & \leq & \begin{cases}
\begin{array}{ll}
\left(1-\frac{1}{\Delta_{H}^{2}}\right)^{d(i,j)} & \textnormal{if }i\textnormal{ is on the boundary of }\Theta_{k}\textnormal{ and }d(i,j)<l_{1}\\
\left(1-\frac{1}{\Delta_{H}^{2}}\right)^{l_{1}-1} & \textnormal{if }i\textnormal{ is on the boundary of }\Theta_{k}\textnormal{ and }d(i,j)=l_{1}\\
0 & \textnormal{\textnormal{otherwise.}}\end{array}\end{cases}\end{eqnarray*}

\end{cor}
\begin{proof}
For each block $\Theta_{k}$ we need to specify a coupling $\Psi_{k}(x,y)$
of the distributions $\Pk(x,\cdot)$ and $\Pk(y,\cdot)$ for each
pair of configurations $(x,y)\in S_{i}$ and each $i\in V$. Trivially
if $i\in\Theta_{k}$ then the set of $H$-colourings with positive
measure in each distribution is the same and the same $H$-colouring
can be chosen for each distribution. The same holds when $i$ is not
on he boundary of $\Theta_{k}$.

Suppose that $i$ is on the boundary of $\Theta_{k}$. Let the other
site on the boundary of $\Theta_{k}$ be coloured $d$ in both $x$
and $y$ and hence $\Pk(x,\cdot)=D_{x_{i},d}^{(l_{1})}$ and $\Pk(y,\cdot)=D_{y_{i},d}^{(l_{1})}$.
We then let $\Psi_{k}(x,y)=\Psi(D_{x_{i},d}^{(l_{1})},D_{y_{i},d}^{(l_{1})})$
which is the coupling constructed in Lemma \ref{lemma:line_coup}
and gives the stated bounds on the disagreement probabilities. 
\end{proof}
\begin{rem*}
It is important to note that, given distinct sites $i$ and $i^{\prime}$
both on the boundary of $\Theta_{k}$, we may use a different coupling
for $\rho_{i,j}^{k}$ and $\rho_{i^{\prime},j}^{k}$. This is the
case since, by definition of $\rho$, the coupling may depend on both
the block and the two initial configurations $x$ and $y$ (which
in turn determine $i$). Since $x$ and $y$ only differ on the colour
assigned to site $i$, the coupling is defined to start from the site
in $\Theta_{k}$ immediately adjacent to $i$, and thus we can use
a different coupling for $\rho_{i,j}^{k}$ and $\rho_{i^{\prime},j}^{k}$. 
\end{rem*}
The following technical lemma is required in the proof of Theorem~\ref{thm:Hcol_2path}.

\begin{lem}
\label{lemma_sum_rank} For any $0\leq p\leq1$ and $j,l\in\Z^{+}$
where $l\geq2j$ \[
p^{j}+p^{l-j+1}\geq p^{j+1}+p^{l-(j+1)+1}.\]

\end{lem}
\begin{proof}
\begin{align*}
p^{j}+p^{l-j+1}-p^{j+1}-p^{l-j} & =p^{j}(1-p)-p^{l-j}(1-p)\\
 & =(p^{j}-p^{l-j})(1-p)\\
 & =p^{j}(1-p^{l-2j})(1-p)\geq0\end{align*}
 since $0\leq p\leq1$ where the last equality uses the fact $l\geq2j$.
\end{proof}
We are now ready to prove Theorem \ref{thm:Hcol_2path}. 

\begin{proof}
[Proof of Theorem \ref{thm:Hcol_2path}] %
\begin{figure}

\caption{A block $\Theta_{k}$ of length $l_{1}$.}

\begin{centering}\label{fig:block-labels}\includegraphics{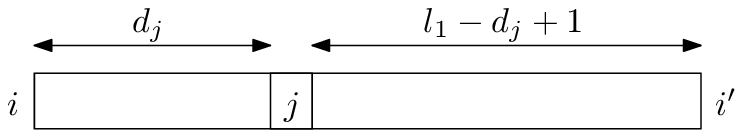}\par\end{centering}
\end{figure}
We will show that $\alpha<1$ and then use Theorem~\ref{thm:main_d}
to obtain the stated bound on the mixing time. Consider some site
$j\in\Theta_{k}$ and let $d_{j}$ denote the number of edges between
$j$ and the nearest site $i\not\in\Theta_{k}$ on the boundary of
$\Theta_{k}$. Then the distance to the other site, $i^{\prime}$,
on the boundary of $\Theta_{k}$ is $l_{1}-d_{j}+1$ as shown in Figure
\ref{fig:block-labels}. Notice that $d_{j}\leq\lceil l_{1}/2\rceil$.
By Corollary \ref{cor:Hcol_rho} we have \[
\rho_{i,j}^{k}\leq\left(1-\frac{1}{\Delta_{H}^{2}}\right)^{d_{j}}\textnormal{ and }\rho_{i^{\prime},j}^{k}\leq\1_{d_{j}\geq2}\left(1-\frac{1}{\Delta_{H}^{2}}\right)^{l_{1}-d_{j}+1}+\1_{d_{j}=1}\left(1-\frac{1}{\Delta_{H}^{2}}\right)^{l_{1}-1}.\]
Now let \[
\alpha_{j,k}=\rho_{i,j}^{k}+\rho_{i^{\prime},j}^{k}\leq\left(1-\frac{1}{\Delta_{H}^{2}}\right)^{d_{j}}+\1_{d_{j}\geq2}\left(1-\frac{1}{\Delta_{H}^{2}}\right)^{l_{1}-d_{j}+1}+\1_{d_{j}=1}\left(1-\frac{1}{\Delta_{H}^{2}}\right)^{l_{1}-1}\]
 be the influence on site $j$. Then \begin{align*}
\alpha & =\max_{k}\max_{j\in\Theta_{k}}\alpha_{j,k}\\
 & \leq\max\left\{ \max_{\left\lceil \frac{l_{1}}{2}\right\rceil \geq d_{j}\geq2}\left\{ \left(1-\frac{1}{\Delta_{H}^{2}}\right)^{d_{j}}+\left(1-\frac{1}{\Delta_{H}^{2}}\right)^{l_{1}-d_{j}+1}\right\} ,\left(1-\frac{1}{\Delta_{H}^{2}}\right)+\left(1-\frac{1}{\Delta_{H}^{2}}\right)^{l_{1}-1}\right\} .\end{align*}
 Since $d_{j}\leq\lceil l_{1}/2\rceil$ the conditions of Lemma \ref{lemma_sum_rank}
are satisfied for $2\leq d_{j}\leq\lceil l_{1}/2\rceil-1$. In particular
taking $d_{j}=\lceil l_{1}/2\rceil-1$, which satisfies the requirements,
gives \[
\left(1-\frac{1}{\Delta_{H}^{2}}\right)^{\lceil l_{1}/2\rceil-1}+\left(1-\frac{1}{\Delta_{H}^{2}}\right)^{l_{1}-\lceil l_{1}/2\rceil+2}\geq\left(1-\frac{1}{\Delta_{H}^{2}}\right)^{\lceil l_{1}/2\rceil}+\left(1-\frac{1}{\Delta_{H}^{2}}\right)^{l_{1}-\lceil l_{1}/2\rceil+1}\]
 and hence \begin{align*}
\max_{\left\lceil \frac{l_{1}}{2}\right\rceil \geq d_{j}\geq2}\left\{ \left(1-\frac{1}{\Delta_{H}^{2}}\right)^{d_{j}}+\left(1-\frac{1}{\Delta_{H}^{2}}\right)^{l_{1}-d_{j}+1}\right\}  & \leq\left(1-\frac{1}{\Delta_{H}^{2}}\right)^{2}+\left(1-\frac{1}{\Delta_{H}^{2}}\right)^{l_{1}-1}\\
 & \leq\left(1-\frac{1}{\Delta_{H}^{2}}\right)+\left(1-\frac{1}{\Delta_{H}^{2}}\right)^{l_{1}-1}\end{align*}
 which gives \begin{align*}
\alpha & \leq\left(1-\frac{1}{\Delta_{H}^{2}}\right)+\left(1-\frac{1}{\Delta_{H}^{2}}\right)^{l_{1}-1}\\
 & =1-\frac{1}{\Delta_{H}^{2}}+\left(1-\frac{1}{\Delta_{H}^{2}}\right)^{\lceil\Delta_{H}^{2}\log(\Delta_{H}^{2}+1)\rceil}\\
 & <1-\frac{1}{\Delta_{H}^{2}}+\frac{1}{\Delta_{H}^{2}+1}\\
 & =1-\frac{1}{\Delta_{H}^{2}(\Delta_{H}^{2}+1)}\end{align*}
 by substituting the definition of $l_{1}$ and using the fact $(1-1/x)^{x}<e^{-1}$
for $x>0$. The statement of the theorem now follows by Theorem \ref{thm:main_d}.
\end{proof}
We now take a moment to show that we are unable to use Theorem~\ref{thm:main_d}
to prove rapid mixing for systematic scan on $H$-colourings of the
$n$-vertex path for any $H$ that does not have a 2-edge path between
all pairs of colours. This motivates the use of path coupling (at
the expense of enforcing a specific scan order) in the subsequent
section.

\begin{observation}\label{obs:noextension} Let $H=(C,E)$ be some
fixed and connected graph in which there is no 2-edge path from $c_{1}$
to $c_{2}$ for some distinct $c_{1},c_{2}\in C$. Then for any set
of $m$ blocks with associated transition matrices $\P1\dots\Pm$
and any coupling $\Psi_{k}(x,y)$ for $1\leq k\leq m$ and $(x,y)\in S_{i}$
we have $\alpha\geq1$ in the unweighted setting. \end{observation}

\begin{proof}
Recall $S_{i}\subseteq\state^{+}\times\state^{+}$ where $\state^{+}$
is the set of all configurations (except when $H$ is bipartite in
which case $\state^{+}$ is one of $\Omega_{1}^{+}$ and $\Omega_{2}^{+}$
as described earlier). Note in particular that any given configuration
in $\state^{+}$ need not be an $H$-colouring of the $n$-vertex
path. Also recall that $\rho_{i,j}^{k}$ is the maximum probability
of disagreement at $j$ when drawing from a coupling starting from
two configurations $(x,y)\in S_{i}$. Let $x$ be any proper $H$-colouring
with $x_{i}=c_{1}$ and $y$ be the configuration with $y_{j}=x_{j}$
for $j\neq i$ and $y_{i}=c_{2}$ (If $H$ is bipartite then $c_{2}$
is from the same colour class of $H$ as $c_{1}$). Note that $y$
is not a proper $H$-colouring as both edges $(y_{i-1},y_{i})\not\in E$
and $(y_{i},y_{i+1})\not\in E$, otherwise the 2-edge paths $(x_{i},x_{i+1}=y_{i+1},y_{i})$
and $(x_{i},x_{i-1}=y_{i-1},y_{i})$ would exist in $H$. However,
$x$ and $y$ are both configurations in $\state^{+}$ and they only
differ at the colour of site $i$ so $(x,y)$ is a valid pair in $S_{i}$.

Now assume that $\alpha<1$. Fix some block $\Theta_{k}=\{ i+1,\dots,i+l\}$
of length $l$ and let $\Pk$ be the transition matrix associated
with $\Theta_{k}$. Also let $\Psi_{k}(x,y)$ be any coupling of $\Pk(x,\cdot)$
and $\Pk(y,\cdot)$. Since $\alpha<1$ it must hold that $\rho_{i,j}^{k}<1$
for each $j\in\Theta_{k}$. In particular $\rho_{i,i+1}^{k}=\Prob_{(x^{\prime},y^{\prime})\in\Psi_{k}(x,y)}(x_{i+1}^{\prime}\neq y_{i+1}^{\prime})<1$
and so (letting $adj(c)$ denote the set of colours adjacent to $c$
in $H$) the set $adj(c_{1})\cap adj(c_{2})$ must be non-empty since
there is a positive probability of assigning the same colour to site
$i+1$ in both distributions. However take any $d\in adj(c_{1})\cap adj(c_{2})$,
then $(c_{1},d,c_{2})$ is a 2-edge path from $c_{1}$ to $c_{2}$
in $H$ contradicting the restriction imposed on $H$ and hence $\alpha\geq1$.
\end{proof}
\begin{rem*}
It remains to be seen if adding weights will allow a proof in the
Dobrushin setting for classes of $H$ not containing 2-edge paths
between all colours. However, this can be done using path coupling
as we will show in section \ref{sec:path_coupling}.
\end{rem*}

\section{$H$-colouring on the path for any $H$\label{sec:path_coupling}}

Recall that $\HscanPC$ is the systematic scan on $\state$ defined
as follows. Let $s=4q+1$, $\beta=\lceil\log(2sq^{s}+1)\rceil q^{s}$
and $l_{2}=2\beta s$. Then $\HscanPC$ is the systematic scan which
performs a heat-bath move on each of the $m_{2}+1=\lfloor2n/l_{2}\rfloor$
blocks in the order $\Theta_{0},\dots,\Theta_{m_{2}}$ where \[
\Theta_{k}=\{ k\beta s+1,\dots,\min((k+2)\beta s,n)\}.\]
Note that the size of $\Theta_{m_{2}}$ is at least $\beta s$ and
that every other block is of size $l_{2}$. We will prove Theorem~\ref{thm:Hcol_line_path}
which bounds the mixing time of $\HscanPC$. Our method of proof will
be path coupling \cite{pathcoupling} and we begin by establishing
some lemmas required to define the coupling we will use in the proof
of Theorem~\ref{thm:Hcol_line_path}. The constructions used in the
following two lemmas are similar to the ones from Lemma 27 in Dyer
et al.~\cite{systematic_scan}.

\begin{lem}
\label{lemma:non-bipartite-path} If $H$ is not bipartite then for
all $c_{1},c_{2}\in C$ there is an $s$-edge path in $H$ from $c_{1}$
to $c_{2}$. 
\end{lem}
\begin{proof}
Let $c\in C$ be some site on an odd-length cycle in $H$ and let
$d_{1}$ be the shortest edge-distance from $c_{1}$ to $c$ and $d_{2}$
the shortest edge-distance from $c$ to $c_{2}$. We construct the
path as follows. Go from $c_{1}$ to $c$ using $d_{1}$ edges. If
$d_{1}+d_{2}$ is even then go around the cycle using an odd number
$q^{\prime}\leq q$ of edges. Go from $c$ to $c_{2}$ in $d_{2}$
edges and observe that the constructed path is of odd length. Also
the length of the path is at most \[
d_{1}+d_{2}+q^{\prime}<3q.\]
Finally go back and forth on the last edge on the path to make the
total length $s$.
\end{proof}
\begin{lem}
\label{lemma:bipartite-path} If $H$ is bipartite with colour classes
$C_{1}$ and $C_{2}$ then for all $c_{1}\in C_{1}$ and $c_{2}\in C_{2}$
there is an $s$-edge path in $H$ from $c_{1}$ to $c_{2}$. 
\end{lem}
\begin{proof}
Go from $c_{1}$ to $c_{2}$ in at most $q-1$ edges and note that
the number of edges is odd. Then go back and forth on the last edge
to make the total path length equal to $s$.
\end{proof}
For completeness we present a proof that $\HscanPC$ is ergodic on
$\state$. 

\begin{lem}
\label{lem:ergodic-H-col-scan}The Markov chain $\HscanPC$ is ergodic
on $\state$.
\end{lem}
\begin{proof}
Let $\HscanPCP$ be the transition matrix of $\HscanPC$. We need
to show that $\HscanPC$ satisfies the following properties 
\begin{itemize}
\item irreducible: $\HscanPCP^{t}(x,y)>0$ for each pair $(x,y)\in\state\times\state$
and some integer $t>0$ 
\item aperiodic: $\gcd\{ t:\HscanPCP^{t}(x,x)>0\}=1$ for each $x\in\state$.
\end{itemize}
In an application of $\HscanPCP$ a heat-bath move is made on each
block in the order $\Theta_{0},\dots,\Theta_{m}$. A heat-bath move
on any block starting from an $H$-colouring has a positive probability
of self-loop which ensures aperiodicity of the chain. To see that
$\HscanPC$ is irreducible consider any pair of $H$-colourings $(x,y)\in\state\times\state$.
We exhibit a sequence of $H$-colourings $x=\sigma^{0},\dots,\sigma^{m_{2}+1}=y$
such that $\sigma_{j}^{k}=\sigma_{j}^{k+1}$ for each $0\leq k\leq m_{2}$
and $j\in V\setminus\Theta_{k}$. Using this sequence we observe that
$\HscanPCP(x,y)>0$ since, for each $0\leq k\leq m_{2}$, performing
a heat-bath move on block $\Theta_{k}$ to $\sigma^{k}\in\state$
results in the $H$-colouring $\sigma^{k+1}\in\state$ with positive
probability. Recall that $\Theta_{k}=\{ k\beta s+1,\dots,\min((k+2)\beta s,n)\}$.
Then let $\sigma^{k}$ be given by \[
\sigma_{i}^{k}=\begin{cases}
y_{i} & \mbox{ if }1\leq i\leq\min((k+2)\beta s-s+1,n)\\
x_{i} & \mbox{ if }(k+2)\beta s+1\leq i\leq n\\
p(i-(k+2)\beta s+s-1) & \mbox{ if }(k+2)\beta s-s+1<i\leq\min((k+2)\beta s,n)\end{cases}\]
where $p(j)$ is the $j$-th in the sequence of colours on the $s$-edge
path in $H$ between $p(0)=y_{(k+2)\beta s-s+1}$ and $p(s)=x_{(k+2)\beta s+1}$
given by Lemmas~\ref{lemma:non-bipartite-path} and~\ref{lemma:bipartite-path}
(since $p(0)$ and $p(s)$ are in opposite colour classes of $H$
in the bipartite case) respectively. 
\end{proof}
The following lemma is an analogue of Lemma 13 in Goldberg et al.~\cite{ssm}.

\begin{lem}
\label{lemma:greedy} For any $c_{1},c_{2},d\in C$ and positive integer
$s^{\prime}\geq s$ such that both $D_{c_{1},d}^{(s^{\prime})}$ and
$D_{c_{2},d}^{(s^{\prime})}$ are non-empty there exists a coupling
$\psi(D_{c_{1},d}^{(s^{\prime})},D_{c_{2},d}^{(s^{\prime})})$ of
$D_{c_{1},d}^{(s^{\prime})}$ and $D_{c_{2},d}^{(s^{\prime})}$ such
that \[
\Prob_{(x^{\prime},y^{\prime})\in\psi(D_{c_{1},d}^{(s^{\prime})},D_{c_{2},d}^{(s^{\prime})})}(x_{v_{s}}^{\prime}\neq y_{v_{s}}^{\prime})\leq1-\frac{1}{q^{s}}.\]

\end{lem}
\begin{proof}
For ease of notation let $D_{1}$ denote $D_{c_{1},d}^{(s^{\prime})}$
and $D_{2}$ denote $D_{c_{2},d}^{(s^{\prime})}$. For $s^{\prime}>s$,
let $n_{k}$ be the number of $H$-colourings on $v_{s+1},\dots,v_{s^{\prime}}$
consistent with $v_{s}$ being assigned colour $k\in C$ and $v_{s^{\prime}}$
adjacent to a site (not in $L$) coloured $d$. If both $s^{\prime}=s$
and $k$ is adjacent to $d$ in $H$ then $n_{k}=1$. If $s^{\prime}=s$
but $k$ is not adjacent to $d$ in $H$ then $n_{k}=0$. The following
definitions are for $i\in\{1,2\}$. Let $l_{i}(k)$ be the number
of $H$-colourings on $v_{1},\dots,v_{s}$ assigning colour $k$ to
site $v_{s}$ and consistent with $v_{1}$ being adjacent to a site
(not in $L$) coloured $c_{i}$. We also let $Z_{i}$ be the set of
$H$-colourings on $L$ with positive measure in $D_{i}$ and $z_{i}$
be the size of this set. Note that $D_{i}$ is the uniform distribution
on $Z_{i}$ so for each $x\in Z_{i}$ $\Pr_{D_{i}}(x)=1/z_{i}$. For
each $k\in C$ let $Z_{i}(k)\subseteq Z_{i}$ be the set of $H$-colourings
with positive measure in $D_{i}$ that assign colour $k$ to site
$v_{s}$ and let $z_{i}(k)$ be the size of this set. Note that $l_{i}(k)n_{k}=z_{i}(k)$
and $\sum_{k}z_{i}(k)=z_{i}$. Let $C_{i}^{*}=\{ k\in C\mid z_{i}(k)>0\}$
be the set of valid colours for $v_{s}$ in $D_{i}$ and let $C^{*}=C_{1}^{*}\cup C_{2}^{*}$.

We define a coupling $\psi$ of $D_{1}$ and $D_{2}$ as follows.
Assume without loss of generality that $z_{1}\geq z_{2}$. We create
the following mutually exclusive subsets of $Z_{i}$. For each $k\in C^{*}$
let $f(k)=\min(z_{1}(k),z_{2}(k))$ and let $F_{1}(k)=\{\sigma^{(k)}(1),\dots,\sigma^{(k)}(f(k))\}\subseteq Z_{1}(k)$
be any subset of $H$-colourings in $Z_{1}$ assigning the colour
$k$ to site $v_{s}$. Also let $F_{2}(k)=\{\tau^{(k)}(1),\dots,\tau^{(k)}(f(k))\}\subseteq Z_{2}(k)$
and observe that $F_{1}(k)$ and $F_{2}(k)$ are of the same size.
We then construct $\psi$ such that for each $k\in C^{*}$ and $j\in\{1,\dots,f(k)\}$
\[
\Prob_{(x^{\prime},y^{\prime})\in\psi}(x^{\prime}=\sigma^{(k)}(j),y^{\prime}=\tau^{(k)}(j))=\frac{1}{z_{1}}.\]
 The rest of the coupling is arbitrary. For example let $R_{i}=Z_{i}\setminus\left(\bigcup_{k\in C^{*}}F_{i}(k)\right)$
be the set of (valid) $H$-colourings not selected in any of the above
subsets of $Z_{i}$ and the size of $R_{i}$ be $r_{i}$, observing
that $r_{1}\geq r_{2}$. Let $R_{1}^{\prime}=\{\sigma(1),\dots,\sigma(r_{2})\}\subseteq R_{1}$
and enumerate $R_{2}$ such that $R_{2}=\{\tau(1),\dots,\tau(r_{2})\}$.
Then for $1\leq j\leq r_{2}$ let \[
\Prob_{(x^{\prime},y^{\prime})\in\psi}(x^{\prime}=\sigma(j),y^{\prime}=\tau(j))=\frac{1}{z_{1}}.\]
 Finish off the coupling by, for each pair $(\sigma\in R_{1}\setminus R_{1}^{\prime},\tau\in Z_{2})$
of $H$-colourings, letting \[
\Prob_{(x^{\prime},y^{\prime})\in\psi}(x^{\prime}=\sigma,y^{\prime}=\tau)=\frac{1}{z_{1}z_{2}}.\]
 From the construction we can verify that the weight of each colouring
$x\in Z_{1}$ in the coupling is $1/z_{1}$ and the weight of each
colouring $y\in Z_{2}$ is \[
\frac{1}{z_{1}}+\frac{z_{1}-z_{2}}{z_{1}z_{2}}=\frac{1}{z_{2}}\]
 since the size of $R_{1}\setminus R_{1}^{\prime}$ is $z_{1}-z_{2}$.
This hence completes the construction of the coupling.

We will require the following bounds on $l_{i}(k)$ for each $k\in C^{*}$
\begin{equation}
1\leq l_{i}(k)\leq q^{s}.\label{eq:firstH}\end{equation}
 There are at most $q$ colours available for each site in the block
and hence at most $q^{s}$ valid $H$-colourings of $v_{1},\dots,v_{s}$
which gives the upper bound. We establish the lower bound by showing
the existence of an $s$-edge path in $H$ from both $c_{1}$ and
$c_{2}$ to any $k\in C^{*}$. Suppose that $H$ is non-bipartite,
then Lemma \ref{lemma:non-bipartite-path} guarantees the existence
of an $s$-edge path in $H$ between any two colours in $H$, satisfying
our requirement.

Now suppose that $H$ is bipartite with colour classes $C_{1}$ and
$C_{2}$. Without loss of generality suppose that $c_{1}\in C_{1}$.
Since both $D_{1}$ and $D_{2}$ are non-empty there exists a $(2s^{\prime}+2)$-edge
path in $H$ from $c_{1}$ to $c_{2}$ (via $d$) so $c_{2}\in C_{1}$.
Let $k\in C$ then $k\in C_{2}$ since there is an $s$-edge path
in $H$ from $c_{1}$ to $k$ and $s$ is odd. Lemma \ref{lemma:bipartite-path}
implies the existence of an $s$-edge path between each $c\in C_{1}$
and each $k\in C_{2}$ which establishes \eqref{eq:firstH}. 

Using \eqref{eq:firstH} to see that $n_{k}\leq f(k)\leq q^{s}n_{k}$
for each $k\in C^{*}$ we have \begin{align*}
\Prob_{(x^{\prime},y^{\prime})\in\psi}(x_{v_{s}}^{\prime}=y_{v_{s}}^{\prime}) & =\sum_{k\in C^{*}}\Prob_{(x^{\prime},y^{\prime})\in\psi}(x_{v_{s}}^{\prime}=y_{v_{s}}^{\prime}=k)\\
 & \geq\sum_{k\in C^{*}}\frac{f(k)}{z_{1}}\\
 & \geq\sum_{k\in C^{*}}\frac{n_{k}}{\sum_{k^{\prime}\in C^{*}}l_{1}(k^{\prime})n_{k}}\\
 & \geq\sum_{k\in C^{*}}\frac{n_{k}}{q^{s}\sum_{k^{\prime}\in C^{*}}n_{k^{\prime}}}\\
 & =\frac{1}{q^{s}}\end{align*}
 which completes the proof.
\end{proof}
\begin{lem}
\label{lemma:small_coup}For any $c_{1},c_{2},d\in C$ and any positive
integer $l^{\prime}\leq l_{2}$ such that both $D_{c_{1},d}^{(l^{\prime})}$
and $D_{c_{2},d}^{(l^{\prime})}$ are non-empty there exists a coupling
$\Psi$ of $D_{c_{1},d}^{(l^{\prime})}$ and $D_{c_{2},d}^{(l^{\prime})}$
in which for $1\leq j\leq l^{\prime}$ \[
\Prob_{(x^{\prime},y^{\prime})\in\Psi(D_{c_{1},d}^{(l^{\prime})},D_{c_{2},d}^{(l^{\prime})})}(x_{v_{j}}^{\prime}\neq y_{v_{j}}^{\prime})\leq\left(1-\frac{1}{q^{s}}\right)^{\left\lfloor \frac{j}{s}\right\rfloor }.\]

\end{lem}
\begin{proof}
We construct a coupling $\Psi(D_{c_{1},d}^{(l')},D_{c_{2},d}^{(l')})$
of $D_{c_{1},d}^{(l^{\prime})}$ and $D_{c_{2},d}^{(l^{\prime})}$
using the following two step process, based on the recursive coupling
in Goldberg et al.~\cite{ssm}. 
\begin{enumerate}
\item If $l^{\prime}<s$ then couple the distributions any valid way which
completes the coupling. Otherwise, couple $D_{c_{1},d}^{(l^{\prime})}(v_{s})$
and $D_{c_{2},d}^{(l^{\prime})}(v_{s})$ greedily to maximise the
probability of assigning the same colour to site $v_{s}$ in both
distributions. Then, independently in each distribution, colour the
sites $v_{1},\dots,v_{s-1}$ consistent with the uniform distribution
on $H$-colourings. Note that it is possible to do this since we obtained
the colour for site $v_{s}$ in each distribution from the induced
distribution on that site. If $l^{\prime}=s$ this completes the coupling. 
\item If the same colour is assigned to $v_{s}$ then the remaining sites
can be coloured the same way in both distributions since the conditional
distributions are the same. Otherwise, for all pairs $(c_{1}^{\prime},c_{2}^{\prime})$
of distinct colours the coupling is completed by recursively constructing
a coupling of $\left[D_{c_{1},d}^{(l^{\prime})}\mid v_{s}=c_{1}^{\prime}\right]=D_{c_{1}^{\prime},d}^{(l^{\prime}-s)}$
and $\left[D_{c_{2},d}^{(l^{\prime})}\mid v_{s}=c_{2}^{\prime}\right]=D_{c_{2}^{\prime},d}^{(l^{\prime}-s)}$. 
\end{enumerate}
This completes the coupling construction and we will prove by strong
induction that for $j\in\{1,\dots,l^{\prime}\}$ \begin{equation}
\Prob_{(x^{\prime},y^{\prime})\in\Psi(D_{c_{1},d}^{(l^{\prime})},D_{c_{2},d}^{(l^{\prime})})}(x_{v_{j}}^{\prime}\neq y_{v_{j}}^{\prime})\leq\left(1-\frac{1}{q^{s}}\right)^{\left\lfloor \frac{j}{s}\right\rfloor }.\label{eq:small_induct}\end{equation}
 Firstly the cases $1\leq j\leq s-1$ are established by observing
that $\lfloor j/s\rfloor=0$ and the probability of disagreement at
any site is at most 1. The case $j=s$ is established in Lemma \ref{lemma:greedy}.
Now for $s<j\leq l^{\prime}$, suppose that \eqref{eq:small_induct}
holds for all positive integers less than $j$. Let $S_{-}=\{ s,2s,\dots\}$
and define the quantities $j_{-}$ and $a_{j}$ by $j_{-}=\max\{ x\in S_{-}\mid x<j\}=a_{j}s$
observing that $1\leq j-j_{-}\leq s$. Now \begin{align*}
\Prob & _{(x^{\prime},y^{\prime})\in\Psi(D_{c_{1},d}^{(l^{\prime})},D_{c_{2},d}^{(l^{\prime})})}(x_{v_{j}}^{\prime}\neq y_{v_{j}}^{\prime})\\
 & =\sum_{c_{1}^{\prime},c_{2}^{\prime}}\Prob_{(x^{\prime},y^{\prime})\in\Psi(D_{c_{1},d}^{(l^{\prime})},D_{c_{2},d}^{(l^{\prime})})}(x_{v_{j_{-}}}^{\prime}=c_{1}^{\prime},y_{v_{j_{-}}}^{\prime}=c_{2}^{\prime})\\
 & \quad\times\Prob_{(x^{\prime},y^{\prime})\in\Psi(D_{c_{1},d}^{(l^{\prime})}\mid v_{j_{-}}=c_{1}^{\prime},D_{c_{2},d}^{(l^{\prime})}\mid v_{j_{-}}=c_{2}^{\prime})}(x_{v_{j}}^{\prime}\neq y_{v_{j}}^{\prime})\\
 & =\sum_{c_{1}^{\prime},c_{2}^{\prime}}\Prob_{(x^{\prime},y^{\prime})\in\Psi(D_{c_{1},d}^{(l^{\prime})},D_{c_{2},d}^{(l^{\prime})})}(x_{v_{j_{-}}}^{\prime}=c_{1}^{\prime},y_{v_{j_{-}}}^{\prime}=c_{2}^{\prime})\\
 & \quad\times\Prob_{(x^{\prime},y^{\prime})\in\Psi(D_{c_{1}^{\prime},d}^{(l^{\prime}-j_{-})},D_{c_{2}^{\prime},d}^{(l^{\prime}-j_{-})})}(x_{v_{j-j_{-}}}^{\prime}\neq y_{v_{j-j_{-}}}^{\prime}).\end{align*}
 Observe that for any pair $(c_{1}^{\prime},c_{2}^{\prime})$ of colours,
if the probabilities of assigning $c_{1}^{\prime}$ to $v_{j_{-}}$
in $D_{c_{1},d}^{(l^{\prime})}$ and $c_{2}^{\prime}$ to $v_{j_{-}}$
in $D_{c_{2},d}^{(l^{\prime})}$ are both non-zero then the distributions
$D_{c_{1}^{\prime},d}^{(l^{\prime}-j_{-})}$ and $D_{c_{2}^{\prime},d}^{(l^{\prime}-j_{-})}$
are both non-empty and hence, using Lemma \ref{lemma:greedy} for
$l^{\prime}-j_{-}\geq s$ and upper-bounding probability of disagreement
by one otherwise, we get \begin{align}
\Prob & _{(x^{\prime},y^{\prime})\in\Psi(D_{c_{1},d}^{(l^{\prime})},D_{c_{2},d}^{(l^{\prime})})}(x_{v_{j}}^{\prime}\neq y_{v_{j}}^{\prime})\nonumber \\
 & \leq\sum_{c_{1}^{\prime},c_{2}^{\prime}}\Prob_{(x^{\prime},y^{\prime})\in\Psi(D_{c_{1},d}^{(l^{\prime})},D_{c_{2},d}^{(l^{\prime})})}(x_{v_{j_{-}}}^{\prime}=c_{1}^{\prime},y_{v_{j_{-}}}^{\prime}=c_{2}^{\prime})\left(\mathbf{1}_{j-j_{-}=s}(1-1/q^{s})+\mathbf{1}_{j-j_{-}\neq s}\right)\nonumber \\
 & \leq\left\{ \begin{array}{ll}
\left(1-\frac{1}{q^{s}}\right)^{\left\lfloor \frac{j_{-}}{s}\right\rfloor +1} & \textnormal{if }j-j_{-}=s\\
\left(1-\frac{1}{q^{s}}\right)^{\left\lfloor \frac{j_{-}}{s}\right\rfloor } & \textnormal{if }j-j_{-}\neq s\end{array}\right.\label{eq:line_floor}\end{align}
 where last inequality is the inductive hypothesis since $j_{-}<j$.

First consider the case $j-j_{-}\neq s$ in which we have $j_{-}+b=j$
for some $1\leq b\leq s-1$. Then \[
\left\lfloor \frac{j_{-}-1}{s}\right\rfloor =\left\lfloor \frac{a_{j}s-1}{s}\right\rfloor =a_{j}-1<a_{j}=\left\lfloor \frac{a_{j}s}{s}\right\rfloor =\left\lfloor \frac{j_{-}}{s}\right\rfloor \]
 and so for $1\leq b\leq s-1$ \[
\left\lfloor \frac{j_{-}+b}{s}\right\rfloor =\left\lfloor \frac{j_{-}}{s}\right\rfloor \]
 which implies that \begin{equation}
\left\lfloor \frac{j_{-}}{s}\right\rfloor =\left\lfloor \frac{j}{s}\right\rfloor .\label{eq:line_less}\end{equation}

Now suppose $j-j_{-}=s$ which substituting for $j_{-}$ gives \begin{equation}
\left\lfloor \frac{j_{-}}{s}\right\rfloor =\left\lfloor \frac{j-s}{s}\right\rfloor =\left\lfloor \frac{j}{s}\right\rfloor -1.\label{eq:line_same}\end{equation}
 Substituting \eqref{eq:line_less} and \eqref{eq:line_same} in \eqref{eq:line_floor}
completes the proof.
\end{proof}
We are now ready to define the coupling of the distributions of configurations
obtained from one complete scan of the Markov chain $\HscanPC$. The
coupling is defined for pairs $(x,y)\in S_{i}$. We will let $(x^{\prime},y^{\prime})$
denote the pair of configurations after one complete scan of $\HscanPC$
starting from $(x,y)$ and let $(x^{k},y^{k})$ be the pair of configurations
obtained by updating blocks $\Theta_{0},\dots,\Theta_{k-1}$ starting
from $(x,y)=(x^{0},y^{0})$. Observe that $(x^{\prime},y^{\prime})$
is obtained by updating block $\Theta_{m_{2}}$ from the pair $(x^{m_{2}},y^{m_{2}})$.

The coupling for updating block $\Theta_{k}$ is defined as follows.
Let $i$ and $i^{\prime}$ be the sites on the boundary of $\Theta_{k}$.
The order of the scan will ensure that at most one of the boundaries
is a disagreement in $(x^{k},y^{k})$, so we only need to define the
coupling for boundaries disagreeing on at most one end of $\Theta_{k}$;
suppose without loss of generality that $x_{i^{\prime}}^{k}=y_{i^{\prime}}^{k}=d$
for some $d\in C$. Firstly, if $x_{i}^{k}=y_{i}^{k}$ then the set
of valid configurations arising from updating $\Theta_{k}$ is the
same in both distributions and we use the identity coupling. 

Otherwise $x_{i}^{k}\neq y_{i}^{k}$. If $H$ is not bipartite then
Lemma~\ref{lemma:non-bipartite-path} implies the existence of a
$(m_{2}+1)$-edge path between both $x_{i}^{k}$ and $d$ and between
$y_{i}^{k}$ and $d$. If $H$ is bipartite then $x_{i}^{k}$ and
$y_{i}^{k}$ are in the same colour class but $d$ is in the opposite
colour class of $H$ since $l_{2}$ is even. Lemma~\ref{lemma:bipartite-path}
implies the existence of a $(m_{2}+1)$-edge path between both $x_{i}^{k}$
and $d$ and between $y_{i}^{k}$ and $d$. Hence both distributions
$D_{x_{i}^{k},d}^{(l_{2})}$ and $D_{y_{i}^{k},d}^{(l_{2})}$ are
non-empty and we obtain $(x^{k+1},y^{k+1})$ from $\Psi(D_{x_{i}^{k},d}^{(l_{2})},D_{y_{i}^{k},d}^{(l_{2})})$
which is the coupling constructed in Lemma~\ref{lemma:small_coup}.
Note that if $k=m_{2}$ (i.e. the block is the last block which may
not be of size $l_{2}$) then both distributions remain (trivially)
non-empty. For ease of reference we state the following corollary
of Lemma~\ref{lemma:small_coup}.

\begin{cor}
\label{cor:Hcol_coupling} For any two sites $v,u\in V$ let $d(v,u)$
denote the edge distance between them. For any block $\Theta_{k}$
let $i$ and $i^{\prime}$ be the sites on the boundary of $\Theta_{k}$
and suppose that $x_{i^{\prime}}^{k}=y_{i^{\prime}}^{k}=d$ for any
$d\in C$. Obtain $(x^{k+1},y^{k+1})$ from the above coupling. Then
for any $j\in\Theta_{k}$ \[
\Pr(x_{j}^{k+1}\neq y_{j}^{k+1})\leq\left\{ \begin{array}{ll}
\left(1-\frac{1}{q^{s}}\right)^{\left\lfloor \frac{d(i,j)}{s}\right\rfloor } & \textnormal{if }x_{i}^{k}\neq y_{i}^{k}\\
0 & \textnormal{otherwise.}\end{array}\right.\]
 
\end{cor}
\begin{lem}
\label{lemma:floor_sum} For any positive integers $s,k,x$ \[
\sum_{j=1}^{sk}\left(1-\frac{1}{x}\right)^{\left\lfloor \frac{j}{s}\right\rfloor }<sx.\]

\end{lem}
\begin{proof}
\[
\sum_{j=1}^{sk}\left(1-\frac{1}{x}\right)^{\left\lfloor \frac{j}{s}\right\rfloor }=(s-1)+s\sum_{j=1}^{k-1}\left(1-\frac{1}{x}\right)^{j}+\left(1-\frac{1}{x}\right)^{k}<s\sum_{j\geq0}\left(1-\frac{1}{x}\right)^{j}<sx.\]

\end{proof}
The following lemma implies Theorem \ref{thm:Hcol_line_path} by Theorem~\ref{thm:path-coupling}
(path coupling).

\begin{lem}
Suppose that $(x,y)\in S_{i}$ and obtain $(x^{\prime},y^{\prime})$
by one complete scan of $\HscanPC$. Then \[
\E{\Ham(x^{\prime},y^{\prime})}<1-\frac{1}{4sq^{s}+2}.\]

\end{lem}
\begin{proof}
First suppose that $i$ is not on the boundary of any block and that
$\Theta_{b}$ is the first block containing $i$. In this case Corollary
\ref{cor:Hcol_coupling} gives us $\Pr(x_{i}^{b+1}\neq y_{i}^{b+1})=0$
and so $\Ham(x^{\prime},y^{\prime})=0.$

Now suppose that $i$ is on the boundary of some block $\Theta_{a}$.
Recall the definition of a block \[
\Theta_{k}=\{ k\beta s+1,\dots,\min(k\beta s+2\beta s,n)\}.\]
 If $i$ is also contained in a block $\Theta_{a^{\prime}}$ with
$a^{\prime}<a$ then Corollary \ref{cor:Hcol_coupling} gives $\Pr(x_{i}^{a^{\prime}+1}\neq y_{i}^{a^{\prime}+1})=0$
and hence $\Ham(x^{\prime},y^{\prime})=0$.

If site $i$ is not updated before $\Theta_{a}$ then $i=(a+2)\beta s+1$
as shown in Figure~\ref{fig:firstPCblock} and the disagreement percolates
through the sites in $\Theta_{a}$ during the update of $\Theta_{a}$.
Using Corollary \ref{cor:Hcol_coupling} we have for $j\in\Theta_{a}$
\begin{equation}
\Pr(x_{j}^{a+1}\neq y_{j}^{a+1})\leq\left(1-\frac{1}{q^{s}}\right)^{\left\lfloor \frac{i-j}{s}\right\rfloor }\label{eq:first_block}\end{equation}
 in particular, the sites in $\Theta_{a}\setminus\Theta_{a+1}=\{ a\beta s+1,\dots(a+1)\beta s\}$
will not get updated again during the scan and hence for $j\in\Theta_{a}\setminus\Theta_{a+1}$
\begin{equation}
\Pr(x_{j}^{\prime}\neq y_{j}^{\prime})\leq\left(1-\frac{1}{q^{s}}\right)^{\left\lfloor \frac{(a+2)\beta s+1-j}{s}\right\rfloor }.\label{eq:final_a}\end{equation}

\begin{figure}

\caption{Site $i$ is on the boundary of $\Theta_{a}$ and is not contained
in any block $\Theta_{a'}$ with $a'<a$.}

\begin{centering}\label{fig:firstPCblock}\includegraphics{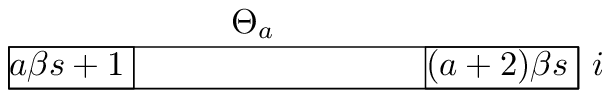}\par\end{centering}
\end{figure}

Now consider the update of any block $\Theta_{k}$ from the pair of
configurations $(x^{k},y^{k})$ where $k>a$. There cannot be a disagreement
at site $(k+2)\beta s+1$ since that site has not been updated (and
it was not the initial disagreement) so the only site on the boundary
of $\Theta_{k}$ that could be a disagreement in $(x^{k},y^{k})$
is $k\beta s$. Hence from Corollary \ref{cor:Hcol_coupling}, for
$j\in\{ k\beta s+1,\dots,\min((k+2)\beta s,n)\}$ \begin{equation}
\Pr(x_{j}^{k+1}\neq y_{j}^{k+1}\mid x_{k\beta s}^{k}\neq y_{k\beta s}^{k})\leq\left(1-\frac{1}{q^{s}}\right)^{\left\lfloor \frac{j-k\beta s}{s}\right\rfloor }.\label{eq:conditional}\end{equation}
 We show by induction on $k$ that for $a+1\leq k\leq m_{2}$ \begin{equation}
\Pr(x_{k\beta s}^{k}\neq y_{k\beta s}^{k})\leq\left(1-\frac{1}{q^{s}}\right)^{\beta(k-a)}.\label{eq:pc_induct}\end{equation}
 The base case, $k=a+1$ follows from \eqref{eq:first_block} since
$j=k\beta s=(a+1)\beta s=a\beta s+\beta s\in\Theta_{a}$. Now suppose
that \eqref{eq:pc_induct} is true for $k-1$. Then \begin{align*}
\Pr(x_{k\beta s}^{k}\neq y_{k\beta s}^{k}) & =\Pr(x_{k\beta s}^{k}\neq y_{k\beta s}^{k}\mid x_{(k-1)\beta s}^{k-1}\neq y_{(k-1)\beta s}^{k-1})\Pr(x_{(k-1)\beta s}^{k-1}\neq y_{(k-1)\beta s}^{k-1})\\
 & \leq\left(1-\frac{1}{q^{s}}\right)^{\left\lfloor \frac{k\beta s-(k-1)\beta s}{s}\right\rfloor }\left(1-\frac{1}{q^{s}}\right)^{\beta(k-a-1)}\\
 & =\left(1-\frac{1}{q^{s}}\right)^{\beta}\left(1-\frac{1}{q^{s}}\right)^{\beta(k-a-1)}\\
 & =\left(1-\frac{1}{q^{s}}\right)^{\beta(k-a)}\end{align*}
 using the inductive hypothesis and \eqref{eq:conditional}.

Now for each site $j\geq(a+1)\beta s+1$, that is site $j$ is updated
at least once following block $\Theta_{a}$, write $j=k_{j}\beta s+b_{j}$
with $1\leq b_{j}\leq\beta s$ where $k_{j}$ denotes is the index
of the block in which $j$ is last updated. \begin{align*}
\Pr(x_{j}^{\prime}\neq y_{j}^{\prime}) & =\Pr(x_{j}^{k_{j}+1}\neq y_{j}^{k_{j}+1})\\
 & \leq\Pr(x_{j}^{k_{j}+1}\neq y_{j}^{k_{j}+1}\mid x_{\beta k_{j}s}^{k_{j}}\neq y_{\beta k_{j}s}^{k_{j}})\Pr(x_{\beta k_{j}s}^{k_{j}}\neq y_{\beta k_{j}s}^{k_{j}}).\end{align*}
 We can then apply \eqref{eq:conditional} to the first component
of the product since $j\in\{ k_{j}\beta s+1,\dots,\min(k_{j}\beta s+2\beta s,n)\}$
and \eqref{eq:pc_induct} to the second since $a+1\leq k_{j}\leq m_{2}$
to get \[
\Pr(x_{j}^{\prime}\neq y_{j}^{\prime})\leq\left(1-\frac{1}{q^{s}}\right)^{\left\lfloor \frac{b_{j}}{s}\right\rfloor }\left(1-\frac{1}{q^{s}}\right)^{\beta(k_{j}-a)}.\]
 Then, using linearity of expectation and \eqref{eq:final_a}, we
have \begin{align*}
\E{\Ham(x^{\prime},y^{\prime})} & =\sum_{j}\Pr(x_{j}^{\prime}\neq y_{j}^{\prime})\\
 & =\sum_{j\in\Theta_{a}\setminus\Theta_{a+1}}\Pr(x_{j}^{\prime}\neq y_{j}^{\prime})+\sum_{j\in\bigcup_{k\geq a+1}\Theta_{k}}\Pr(x_{j}^{\prime}\neq y_{j}^{\prime})\\
 & \leq\sum_{j=as\beta+1}^{(a+1)\beta s}\left(1-\frac{1}{q^{s}}\right)^{\left\lfloor \frac{(a+2)\beta s+1-j}{s}\right\rfloor }+\sum_{k_{j}=a+1}^{m_{2}}\sum_{b_{j}=1}^{\beta s}\left(1-\frac{1}{q^{s}}\right)^{\left\lfloor \frac{b_{j}}{s}\right\rfloor }\left(1-\frac{1}{q^{s}}\right)^{\beta(k_{j}-a)}\\
 & =\sum_{r=1}^{\beta s}\left(1-\frac{1}{q^{s}}\right)^{\left\lfloor \frac{\beta s+r}{s}\right\rfloor }+\sum_{k_{j}=a+1}^{m_{2}}\left(1-\frac{1}{q^{s}}\right)^{\beta(k_{j}-a)}\sum_{b_{j}=1}^{\beta s}\left(1-\frac{1}{q^{s}}\right)^{\left\lfloor \frac{b_{j}}{s}\right\rfloor }\\
 & <\left(1-\frac{1}{q^{s}}\right)^{\beta}\sum_{r=1}^{\beta s}\left(1-\frac{1}{q^{s}}\right)^{\left\lfloor \frac{r}{s}\right\rfloor }+\sum_{t\geq1}\left(\left(1-\frac{1}{q^{s}}\right)^{\beta}\right)^{t}\sum_{b_{j}=1}^{\beta s}\left(1-\frac{1}{q^{s}}\right)^{\left\lfloor \frac{b_{j}}{s}\right\rfloor }\\
 & <\left(1-\frac{1}{q^{s}}\right)^{\beta}sq^{s}+\frac{\left(1-\frac{1}{q^{s}}\right)^{\beta}sq^{s}}{1-\left(1-\frac{1}{q^{s}}\right)^{\beta}}\end{align*}
 where the last inequality uses Lemma \ref{lemma:floor_sum} and the
sum of a geometric progression. Substituting the definition of $\beta$
and using the fact $(1-1/x)^{x}<e^{-1}$ for $x>0$ we get \begin{align*}
\E{\Ham(x^{\prime},y^{\prime})} & <\left(1-\frac{1}{q^{s}}\right)^{\lceil\log(2sq^{s}+1)\rceil q^{s}}sq^{s}+\frac{\left(1-\frac{1}{q^{s}}\right)^{\lceil\log(2sq^{s}+1)\rceil q^{s}}sq^{s}}{1-\left(1-\frac{1}{q^{s}}\right)^{\lceil\log(2sq^{s}+1)\rceil q^{s}}}\\
 & <\frac{sq^{s}}{e^{\lceil\log(2sq^{s}+1)\rceil}}+\frac{sq^{s}}{e^{\lceil\log(2sq^{s}+1)\rceil}(1-e^{-\lceil\log(2sq^{s}+1)\rceil})}\\
 & =\frac{sq^{s}}{e^{\lceil\log(2sq^{s}+1)\rceil}}+\frac{sq^{s}}{e^{\lceil\log(2sq^{s}+1)\rceil}-1}\\
 & \leq\frac{sq^{s}}{2sq^{s}+1}+\frac{sq^{s}}{2sq^{s}}\\
 & =1-\frac{1}{4sq^{s}+2}\end{align*}
 which completes the proof.
\end{proof}

\section{$H$-colouring using a random update Markov chain\label{sec:random-update}}

Recall that the random update Markov chain $\Hrandom$ on $\state$
is defined as follows. We again let $s=4q+1$ and we define $\gamma=2q^{s}+1$.
We then define a set of $n+s\gamma-1$ blocks of size at most $s\gamma$
as follows. \[
\Theta_{k}=\begin{cases}
\{ k,\dots,\min(k+s\gamma-1,n)\} & \textnormal{when }k\in\{1,\dots,n\}\\
\{1,\dots,n+s\gamma-k\} & \textnormal{when }k\in\{ n+1,\dots,n+s\gamma-1\}\end{cases}\]
By construction of the set of blocks each site is adjacent to \emph{at
most} two blocks and furthermore each site is contained in \emph{exactly}
$s\gamma$ blocks. One step of $\Hrandom$ consists of selecting a
block uniformly at random and performing a heat-bath update on it.
We will prove (using path coupling) Theorem~\ref{thm:random_Hcolor}
namely that $\Hrandom$ mixes in $O(n\log n)$ updates for any $H$.

We begin by defining the required coupling. For a pair of configurations
$(x,y)\in S_{i}$ we obtain the pair $(x^{\prime},y^{\prime})$ by
one step of $\Hrandom$. That is we select a block uniformly at random
and perform a heat bath move on that block. We can again use Lemma~\ref{lemma:small_coup}
from Section~\ref{sec:path_coupling} to construct the required coupling
for updating block $\Theta_{k}$ since the definition of $s$ is the
same in both Markov chains. If $i$ is not on the boundary of $\Theta_{k}$
then the sets of valid $H$-colourings of $\Theta_{k}$ are the same
in both distributions and we use the identity coupling. If $i$ is
on the boundary of $\Theta_{k}$ then we let the other site on the
boundary be coloured $d$ in both $x$ and $y$. We then obtain $(x^{\prime},y^{\prime})$
from $\Psi(D_{x_{i},d}^{(s\gamma)},D_{y_{i},d}^{(s\gamma)})$ which
is the coupling constructed in Lemma~\ref{lemma:small_coup}. The
disagreement probabilities are summarised in the following corollary
(of Lemma~\ref{lemma:small_coup}).

\begin{cor}
\label{cor:random_coupling} For any two sites $v,u\in V$ let $d(v,u)$
denote the edge distance between them. Suppose that a block $\Theta_{k}$
has been selected to be updated. For any pair $(x,y)\in S_{i}$ obtain
$(x^{\prime},y^{\prime})$ from the above coupling. Then for any $j\in\Theta_{k}$
\[
\Pr(x_{j}^{\prime}\neq y_{j}^{\prime})\leq\left\{ \begin{array}{ll}
\left(1-\frac{1}{q^{s}}\right)^{\left\lfloor \frac{d(i,j)}{s}\right\rfloor } & \textnormal{if }i\textnormal{ is on the boundary of }\Theta_{k}\\
0 & \textnormal{otherwise.}\end{array}\right.\]

\end{cor}
The following lemma implies Theorem~\ref{thm:random_Hcolor} by Theorem~\ref{thm:path-coupling}
(path coupling).

\begin{lem}
Suppose that $(x,y)\in S_{i}$ and obtain $(x^{\prime},y^{\prime})$
by one step of $\Hrandom$. Then \[
\E{\Ham(x^{\prime},y^{\prime})}<1-\frac{s}{n+2sq^{s}+s-1}.\]

\end{lem}
\begin{proof}
There are $s\gamma$ blocks containing site $i$ and if such a block
is selected then $\Ham(x^{\prime},y^{\prime})=0$. There are at most
2 blocks adjacent to site $i$ and if such a block is selected then
the discrepancy percolates in the block according to the probabilities
stated in Corollary~\ref{cor:random_coupling}. This leaves $n+s\gamma-1-s\gamma-2=n-3$
blocks that leave the Hamming distance unchanged. Hence, using Lemma~\ref{lemma:floor_sum},
we have\begin{eqnarray*}
\E{\Ham(x^{\prime},y^{\prime})} & \leq & \frac{2}{n+s\gamma-1}\left(1+\sum_{j=1}^{\gamma s}\left(1-\frac{1}{q^{s}}\right)^{\left\lfloor \frac{j}{s}\right\rfloor }\right)+\frac{n-3}{n+s\gamma-1}\\
 & < & \frac{n-1}{n+s\gamma-1}+\frac{2sq^{s}}{n+s\gamma-1}\\
 & = & \frac{2sq^{s}+n-1}{2sq^{s}+n-1+s}=1-\frac{s}{2sq^{s}+n-1+s}\end{eqnarray*}
by substituting the definition of $\gamma.$
\end{proof}

\section*{Acknowledgments}

I am grateful to Leslie Goldberg for several useful discussions regarding
technical issues and for providing detailed and helpful comments on
a draft of this article.\bibliographystyle{plain}
\bibliography{../references}

\begin{thebibliography}{10}

\bibitem{grid_ach}
Dimitris Achlioptas, Mike Molloy, Cristopher Moore, and Frank {Van~Bussel}.
\newblock Sampling grid colourings with fewer colours.
\newblock In {\em Proc.~of the 6th Latin American Symposium on Theoretical
  Informatics (LATIN'04)}, pages 80--89, Buenos Aires, Argentina, 2004.

\bibitem{aldous_walks}
David~J Aldous.
\newblock Random walks on finite groups and rapidly mixing markov chains.
\newblock In {\em S\'eminaire de probabilit\'es XVII}, pages 243--297.
  Springer-Verlag, 1983.

\bibitem{amit}
Yali Amit.
\newblock Convergence properties of the {Gibbs} sampler for pertubations of
  gaussians.
\newblock {\em The Annals of Statistics}, 24(1):122--140, 1996.

\bibitem{metric}
Magnus Bordewich, Martin Dyer, and Marek Karpinski.
\newblock Stopping times, metrics and approximate counting.
\newblock In Michele Bugliesi, Bart Preneel, Vladimiro Sassone, and Ingo
  Wegener, editors, {\em ICALP}, volume 4051 of {\em Lecture Notes in Computer
  Science}, pages 108--119. Springer, 2006.

\bibitem{pathcoupling}
Russ Bubley and Martin Dyer.
\newblock Path coupling: a technique for proving rapid mixing in {Markov}
  chains.
\newblock In {\em 38th Annual Symposium on Foundations of Computer Science},
  pages 223--231, 1997.

\bibitem{beach}
Robert Burton and Jeffrey Steif.
\newblock Nonuniqueness of measures of maximal entropy for subshifts of finite
  type.
\newblock {\em Ergodic Theory and Dynamical Systems}, 14(2):213--236, 1994.

\bibitem{cooper}
Colin Cooper, Martin Dyer, and Alan Frieze.
\newblock On {Markov} chains for randomly {$H$}-colouring a graph.
\newblock {\em Journal of Algorithms}, 39(1):117--134, 2001.

\bibitem{MCMCDiagnostics}
Mary~Kathryn Cowles and Bradlet~P. Carlin.
\newblock Markov chain {Monte Carlo} convergence diagnostics: A comparative
  review.
\newblock {\em Journal of The American Statistical Association},
  91(434):883--904, 1996.

\bibitem{diaconis-scan}
Persi Diaconis and Arun Ram.
\newblock Analysis of systematic scan {Metropolis} algorithms using
  {Iwahoti-Hecke} algebra techniques.
\newblock {\em Michigan Mathematical Journal}, 48:157--190, 2000.

\bibitem{ind_set_sparse}
Martin Dyer, Alan Frieze, and Mark Jerrum.
\newblock On counting independent sets in sparse graphs.
\newblock {\em SIAM Journal Computing}, 31(5):1527--1541, 2002.

\bibitem{relative}
Martin Dyer, Leslie~Ann Goldberg, Catherine Greenhill, and Mark Jerrum.
\newblock On the relative complexity of approximate counting problems.
\newblock {\em Algorithmica}, 38(3):471--500, 2003.

\bibitem{H-col-sample}
Martin Dyer, Leslie~Ann Goldberg, and Mark Jerrum.
\newblock Counting and sampling {$H$}-colourings.
\newblock {\em Information and Computation}, 189:1--16, 2004.

\bibitem{dobrushin_scan}
Martin Dyer, Leslie~Ann Goldberg, and Mark Jerrum.
\newblock Dobrushin conditions and systematic scan.
\newblock In Josep D\'{\i}az, Klaus Jansen, Jos{\'e} D.~P. Rolim, and Uri
  Zwick, editors, {\em APPROX-RANDOM}, volume 4110 of {\em Lecture Notes in
  Computer Science}, pages 327--338. Springer, 2006.

\bibitem{systematic_scan}
Martin Dyer, Leslie~Ann Goldberg, and Mark Jerrum.
\newblock Systematic scan and sampling colourings.
\newblock {\em Annals of Applied Probability}, 16(1):185--230, 2006.

\bibitem{H-col-count}
Martin Dyer and Catherine~S. Greenhill.
\newblock The complexity of counting graph homomorphisms.
\newblock {\em Random Structures and Algorithms}, 17:260--289, 2000.

\bibitem{ind_set}
Martin Dyer and Catherine~S. Greenhill.
\newblock On {Markov} chains for independent sets.
\newblock {\em J. Algorithms}, 35(1):17--49, 2000.

\bibitem{H-col-bounded}
Anna Galluccio, Pavol Hell, and Jaroslav~Ne\v set\v ril.
\newblock The complexity of {$H$}-colouring of bounded degree graphs.
\newblock {\em Discrete Mathematics}, 222:101--109, 2000.

\bibitem{kelk}
Leslie~Ann Goldberg, Steven Kelk, and Mike Paterson.
\newblock The complexity of choosing an {$H$}-colouring (nearly) uniformly at
  random.
\newblock {\em SICOMP}, 33(2):416--432, 2004.

\bibitem{ssm}
Leslie~Ann Goldberg, Russ Martin, and Mike Paterson.
\newblock Strong spatial mixing for lattice graphs with fewer colours.
\newblock {\em SICOMP}, 35(2):486--517, 2005.

\bibitem{H-col-NP}
Pavol Hell and Jaroslav~Ne\v set\v ril.
\newblock On the complexity of {$H$}-colouring.
\newblock {\em Journal of Combinatorial Theory, Series B}, 48:92--110, 1990.

\bibitem{jerrum_simple}
Mark Jerrum.
\newblock A very simple algorithm for estimating the number of $k$-colourings
  of a low-degree graph.
\newblock {\em Random Structures and Algorithms}, 1995.

\bibitem{kenyon_tree}
Claire Kenyon, Elchanan Mossel, and Yuval Peres.
\newblock Glauber dynamics on trees and hyperbolic graphs.
\newblock In {\em Proc. 42nd Annual IEEE Symposium on Foundations of Computer
  Science}, pages 568--578, 2001.

\bibitem{luby-vigoda-indset-improve}
Michael Luby and Eric Vigoda.
\newblock Fast convergence of the {G}lauber dynamics for sampling independent
  sets: {P}art {I}.
\newblock {\em Random Structures and Algorithms}, 15(3--4):229--241, 1999.

\bibitem{tree}
Fabio Martinelli, Alistair Sinclair, and Dror Weitz.
\newblock Glauber dynamics on trees: Boundary conditions and mixing time.
\newblock {\em Communications in Mathematical Physics}, 250(2):301--334, 2004.

\bibitem{dobrushin_ee}
Kasper Pedersen.
\newblock Dobrushin conditions for systematic scan with block dynamics
  (extended abstract).
\newblock To appear in MFCS, 2007.

\bibitem{salas-sokal}
Jesus Salas and Alan~D Sokal.
\newblock Absence of phase transition for antiferromagnetic potts models via
  the dobrushin uniqueness theorem.
\newblock {\em Journal of Statistical Physics}, pages 551--579, 1997.

\bibitem{vigoda}
Eric Vigoda.
\newblock Improved bounds for sampling colourings.
\newblock {\em J. Math. Phys}, 2000.

\bibitem{dror_combinatorial}
Dror Weitz.
\newblock Combinatorial criteria for uniqueness of {Gibbs} measures.
\newblock {\em Random Structures and Algorithms}, 27(4):445--475, 2005.

\bibitem{dror_ind_set}
Dror Weitz.
\newblock Counting independent sets up to the tree threshold.
\newblock In {\em STOC}, pages 140--149, 2006.

\bibitem{WR}
Benjamin Widom and John~S. Rowlinson.
\newblock New model for the study of liquid-vapour phase transition.
\newblock {\em The Journal of Chemical Physics}, 52(4):1670--1684, 1970.

\end{thebibliography}

\end{document}